\documentclass{amsart}[14pt]

\usepackage{amsmath, amstext,  amssymb, amsopn, amsthm, amsbsy, amsgen, amscd, amsfonts, fouriernc}

\newcommand{\A}{\mathbb A}
\newcommand{\B}{\mathbb B}
\newcommand{\R}{\mathbb R}
\newcommand{\C}{\mathbb C}

\newcommand{\HP}{\mathbb H}

\newcommand{\N}{\mathbb N}
\newcommand{\PR}{\mathbb P}
\newcommand{\Q}{\mathbb Q}
\newcommand{\Z}{\mathbb Z}

\newcommand{\T}{\mathbb T}

\newcommand{\F}{\mathbb F}

\newcommand{\SO}{\hat{\mathbb S}}
\newcommand{\HSO}{\hat{\mathfrak{S}}}

\newtheorem{theo}{Theorem}
\newtheorem{lemm}{Lemma}
\newtheorem{prop}{Proposition}
\newtheorem{coro}{Corollary}

\newtheorem*{stheo2}{Theorem}
\newtheorem*{wigtheo}{Wigner's Theorem}

\theoremstyle{definition}
\newtheorem{defi}{Definition}

\theoremstyle{remark}

\newtheorem{note}{Note}

\newtheorem{exam}{Example}

\title[Geometric Galois Theory]{Geometric Galois Theory, Nonlinear Number Fields
and a Galois Group Interpretation of the Idele Class Group \\
{\small {\sc Revised Version}}}
\author
{T.M. Gendron \& A. Verjovsky}
\address{Instituto de Matematicas -- Unidad Cuernavaca, Universidad
Autonoma Nacional de M\'{e}xico, Av. Universidad S/N, C.P. 62210
Cuernavaca, Morelos, MEXICO}
\date{12 April, 2010}
\dedicatory{In memory of Egidio Barrera.}
\begin{document}
\maketitle
\begin{abstract} This paper concerns the description of
holomorphic extensions of algebraic number fields. We
define a hyperbolized adele class group $\HSO_{K}$ for every algebraic number field $K$
and consider the Hardy space ${\sf H}_{\bullet}[K]$
of graded-holomorphic
functions on $\HSO_{K}$.  We show that the subspace 
${\sf N}_{\bullet}[K]$ of the projectivization
$\PR{\sf H}_{\bullet}[K]$ defined by
the functions of non-zero trace possesses
two partially-defined operations $\oplus$ and $\otimes$, with respect
to which there is canonical monomorphism
of $K$ into ${\sf N}_{\bullet}[K]$.  We call ${\sf N}_{\bullet}[K]$ a {\it nonlinear
field extension} of $K$.
 We define Galois groups
for nonlinear fields and show that
${\sf Gal}({\sf N}_{\bullet}[L]/{\sf N}_{\bullet}[K])$ $\cong$
${\sf Gal}(L/K)$ if $L/K$ is Galois.  If $\Q^{\rm ab}$ denotes the maximal abelian
extension of $\Q$, ${\sf C}_{\Q}$  the idele class
group and $\bar{\sf N}[\Q^{\rm ab}]=\PR{\sf H}_{\bullet}[K]$ is the full projectivization 
then we show that there
exist embeddings of ${\sf C}_{\Q}$ into
${\sf Gal}_{\oplus}(\bar{\sf N}_{\bullet}[\Q^{\rm ab}]/\Q)$
and 
${\sf Gal}_{\otimes}(\bar{\sf N}_{\bullet}[\Q^{\rm ab}]/\Q)$, the ``Galois
groups'' of automorphisms preserving $\oplus$ resp.\ $\otimes$ only. 
\end{abstract}
\tableofcontents

\section{Introduction}

Since the introduction of a global understanding of number theory (Gauss, Galois) and geometry
(Riemann), the idea that the two subjects exist in parallel duality  
has exercised a tremendous pull on mathematical thought.
An indication of this conjectural relationship can be found in 
the equivalence between (coverings of) Riemann surfaces and (extensions of)
their fields of meromorphic functions, based on which, one
can formulate the following meta-principle:

\begin{quote}  {\textit {\small To every algebraic number field $K/\Q$, 
there exists a ``Riemann surface'' $\Sigma_{K}$ for which
\[   ``{\sf Mer}(\Sigma_{K})\;\; \cong \;\; K\text{''} .    \]
If $L/K$ is a Galois extension, then $\Sigma_{L}$ 
is a Galois covering of $\Sigma_{K}$ and
\[ ``{\sf Deck}(\Sigma_{L}/\Sigma_{K})\;\;\cong\;\; {\sf Gal}(L/K )\text{''}.\]}}
\end{quote}

If we were interested in an extension $E$ over the function field $\F_{p^{n}}[X]$, then
this principle is, suitably interpreted, correct:
there exists a curve over $\F_{p^{n}}$ whose function field is $E$.  This observation
was used by Weil to prove the Riemann hypothesis for function fields over finite fields \cite{We1}.

In the case of a number field $K$, 
the principle as stated
above must be modified, as any reasonable
notion of ``field of meromorphic functions'' must contain
the field of constants $\C$, which  cannot be a subfield of
$K$.  Thus we might instead
ask that $K$ generate ${\sf Mer}(\Sigma_{K})$ over $\C$.  
The additional freedom afforded by the use of
$\C$-coefficients motivates
a second meta-principle:

\begin{quote}
{\textit {\small  If
$K^{\rm ab}$ is the maximal abelian extension
of $K$ and ${\sf C}_{K}$ is the idele class group of $K$, then 
there is a monomorphism 
\[ ``{\sf C}_{K}\;\hookrightarrow\; {\sf Gal}\big({\sf Mer}(\Sigma_{K^{\rm ab}})/K\big)\text{''}.\]
}}
\end{quote}

This second meta-principle -- which is true for function fields over finite fields 
upon using rational function fields in place of meromorphic function fields --  
has an important place in Weil's approach to the classical
Riemann hypothesis, as the following often quoted passage reveals
\cite{We2}:

\begin{quotation}  {\small ``La recherche d'une interpr\'{e}tation pour
${\sf C}_{K}$ si $K$ est un corps de nombres, analogue en
quelque mani\`{e}re \`{a} l'interpr\'{e}tation par un groupe de Galois
quand $K$ est un corps de fonctions, me semble constituer l'un
des probl\`{e}mes fondamentaux de la th\'{e}orie des nombres \`{a} l'heure
actuelle; il se peut qu'une telle interpr\'{e}tation renferme la
clef de l'hypoth\`{e}se Riemann...''}
\end{quotation}

Weil's speculation regarding a Galois group interpretation of ${\sf C}_{K}$
has inspired new
approaches to the Riemann hypothesis, {\em e.g.}\ especially that of Alain
Connes \cite{Co}.  In this paper, we shall give a certain expression to
these principles through a hyperbolized version of the 
adele class group of $K$.

Let us consider first a field $K$ of finite degree $d$ over $\Q$.
To $K$ we may associate the adele class group $\SO_{K}=\A_{K}/(K,+)$, 
a $d$-dimensional solenoid 
whose leaves are dense and isomorphic to a product of the form
 \[ K_{\infty}=\R^{r}\times\C^{s}\] where $K$ has $r$ real places and $2s$ complex places. 
From $\SO_{K}$ we construct a hyperbolization $\HSO_{K}$
whose leaves are polydisks isomorphic to $(\HP^{2})^{d}$ and whose distinguished
boundary is $\SO_{K}$.  See \S\S \ref{SolenoidChap}, \ref{AdClassIChap}  and \ref{hyperbolizations} for more details.

For a number field $\mathcal{K}/\Q$ of infinite degree,
such as the maximal abelian extension $K^{\rm ab}$, the notion of adele class group
has not yet, to our knowledge, been defined. 
In order to redeem the desired properties found in the
adele class group of a finite field extension, it is necessary to 
consider a {\em pair} of adele class groups 
that work in tandem.  Representing $\mathcal{K}=\underset{\longrightarrow}{\lim}\, K_{\lambda}$, where the 
$K_{\lambda}$ are finite degree extensions of $\Q$, 
the (ordinary) adele class group $\SO_{\mathcal{K}}$ is formed from
$\underset{\longrightarrow}{\lim} \,\SO_{K_{\lambda}}$
by completing its canonical
leaf-wise euclidean metric, see \S 4.  Since
$\SO_{\mathcal{K}}$ is not locally-compact, we 
consider in \S 5 a compactification $\hat{\SO}_{\mathcal{K}}$ 
called the {\em proto adele class group} of $\mathcal{K}$, arising 
as the inverse limit of the trace maps.
We use the hilbertian $\SO_{\mathcal{K}}$ to
define the hyperbolization $\HSO_{\mathcal{K}}$, and the compact
$\hat{\SO}_{\mathcal{K}}$ to provide Fourier analysis.  See \S\S \ref{AdClassIIChap}, \ref{ProtoAdClassChap}.

The origin of the notion of a nonlinear field
comes from an enhanced understanding of
the character group ${\sf Char}\big(\SO_{K}\big)$, described in \S \ref{CharFieldChap}.  
The character group 
possesses an additional operation making it a field, 
and there is a canonical isomorphism
\begin{equation}\label{dual} {\sf Char}\big(\SO_{K}\big)\;\; \cong\;\; K .
\end{equation}
If we denote by $\T_{K}=K_{\infty}/O_{K}$ the Minkowski torus, the above isomorphism
identifies 
\[  {\sf Char}\big(\T_{K}\big)\;\; \cong \;\;\mathfrak{d}_{K}^{-1} \supset O_{K} \]
where $\mathfrak{d}_{K}^{-1}$ is the inverse different.  Thus ${\sf Char}\big(\T_{K}\big)$ has the enhanced
structure of an $O_{K}$-module canonically extending the ring $O_{K}$.
The same
is true for an infinite field extension $\mathcal{K}$ that is Galois over $\Q$ if we use the proto adele class group $\hat{\SO}_{\mathcal{K}}$.

Let $f\in L^{2}(\SO_{K}, \, \C )$ = the Hilbert space of square integrable
complex-valued functions with respect to normalized Haar measure.  
By Fourier theory, $f$ has the development
\[f=\sum a_{\alpha}\phi_{\alpha}\]
where $\alpha\in K$, $a_{\alpha}\in\C$ and $\phi_{\alpha}\in  {\sf Char}\big(\SO_{K}\big)$, and so
the isomorphism (\ref{dual}) defines an inclusion
$K\hookrightarrow L^{2}(\SO_{K}, \, \C )$.
Cauchy (point-wise) multiplication of functions, when defined, is denoted
$f\oplus g$ since it restricts to $+_{K}$ in $K$.  
The operation corresponding to $\times_{K}$
is the Dirichlet product $f\otimes g$, and through it
$L^{2}(\SO_{K}, \, \C )$ acquires
the structure of a {\it partial double group algebra}, where partial refers to the fact
that the two operations are only defined when square integrability is conserved.  The departure from ordinary
field theory begins with the observation
that Dirichlet multiplication {\em does not} distribute over Cauchy
multiplication, or to put it differently, the extension
of multiplication in $K$ to $L^{2}(\SO_{K},\, \C )$ no
longer defines a Cauchy bilinear operation.  See \S \ref{GradedChap}.

An interpretation of $L^{2}(\SO_{K},\, \C )$ as boundary values of
holomorphic functions on $\HSO_{K}$ comes 
through the introduction of graded holomorphicity, treating, in the style of conformal field theory,  each notion of holomorphicity
on the same footing as the classical one.
Suppose first that $K$ is totally real. Denote by $\Theta_{K}=\{ -, +\}^{d}$ the sign group.
To each $\theta\in\Theta_{K}$ we associate the Hardy space ${\sf H}_{\theta}[K]$ of $\theta$-holomorphic functions on $\HSO_{K}$, and every $f\in L^{2}(\SO_{K},\, \C )$ determines
a $2^{d}+1$ tuple $(F_{\theta}; F_{0})$ consisting of $\theta$-holomorphic components and the constant term $F_{0}=a_{0}$.  In particular,
there is an isomorphism of Hilbert spaces
\[  {\sf H}_{\bullet}[K] := \bigoplus_{\theta}  {\sf H}_{\theta}[K]\oplus\C \cong  L^{2}(\SO_{K},\, \C ). \]
The graded Hilbert space ${\sf H}_{\bullet}[K]$ inherits
the partially-defined operations of $\oplus$ and $\otimes$, where the Dirichlet product has a homogeneous decomposition with respect to the grading of ${\sf H}_{\bullet}[K]$.  See \S \ref{GradedSectReal}.

When $K$ is totally complex we consider, for each complex place, a pair of order four signings:
the {\it singular complex sign group} $\Theta^{\C}=\{ \sqrt{-},-,-  \sqrt{-},+ \}$, which signs points on $(\R\cup i\R)-0$ according to the axial component
that they belong to, 
and the  {\it  complex sign set} $\Omega=\{ \sqrt{-}\epsilon ,-\epsilon ,-  \sqrt{-}\epsilon ,+\epsilon \}$
which signs points in $\C -(\R\cup i\R)$ according to the quadrant they belong to.  If we denote by
$\mho=\Theta^{\C}\cup \Omega$ then there is 
a subset $\mho_{K}\subset \mho^{s}$ with respect to which we may define, for each $\vartheta\in\mho_{K}$, a Hardy space ${\sf H}_{\vartheta}[K]$ of $\vartheta$-holomorphic functions on
$\HSO_{K}$.  We obtain in analogy with the real case an isomorphism
\[  {\sf H}_{\bullet}[K] := \bigoplus_{\vartheta}  {\sf H}_{\vartheta}[K]\oplus\C \cong  L^{2}(\SO_{K},\, \C ). \]
See \S \ref{GradedSectComplex}.
In the hybrid case where $K$ has both real and complex planes we obtain a direct sum
${\sf H}_{\bullet}[K]= {\sf H}_{\bullet}^{\R} [K]\oplus {\sf H}^{\C}_{\bullet}[K]$ consisting of the
$\theta$-holomorphic and $\vartheta$-holomorphic parts.

Since the vector space structure of ${\sf H}_{\bullet}[K]$
does not descend to $K$, it is natural to discard it by projectivizing.
After removing the functions whose trace is defined and equal to zero,
we obtain an infinite-dimensional affine subspace
\[ {\sf N}_{\bullet}[K]\subset \PR{\sf H}_{\bullet}[K] \]
endowed with a partial double semigroup structure
induced from $\oplus$ and 
and $\otimes$, satisfying the following properties:
\begin{enumerate}
\item[1.] Let $\C [K]$ be the field algebra (double group algebra) associated to $K$ and let ${\sf N}^{0}_{\bullet}[K]\subset \PR \C[K]$ denote the sub double semigroup of 
elements of non zero trace, graded according to the same scheme described above.  Then
there is a graded monomorphism 
${\sf N}^{0}_{\bullet}[K]\hookrightarrow 
{\sf N}_{\bullet}[K]$ with dense image.
\item[2.]  The identity ${\sf id}_{\oplus}$ is a universal annihilator
for the product $\otimes$: for all $f\in {\sf N}_{\bullet}[K]$ for which $f\otimes {\sf id}_{\otimes}$
is defined, $f\otimes {\sf id}_{\otimes}={\sf id}_{\otimes}$.
\end{enumerate}
We call a topological partial double semi-group satisfying these 
properties an {\em (abstract) nonlinear number field} over $K$.  
The ring of integers $O_{K}$ of $K$ generates
in turn the {\em nonlinear ring of integers}
${\sf N}_{\bullet}[O_{K}]\;\subset\;
{\sf N}_{\bullet}[K]$, a nonlinear extension of $O_{K}$.  On
a dense subset of ${\sf N}_{\bullet}[K]$, every element
may be regarded as a Dirichlet quotient of elements of 
${\sf N}_{\bullet}[O_{K}]$.  See \S \ref{NonlinearChap}.

An automorphism of the nonlinear number field ${\sf N}_{\bullet}[K]$
is defined to be the restriction of a graded Fubini-Study isometry
preserving the operations $\oplus$ and $\otimes$. Given $K$ an algebraic number field,
denote by
\[ {\sf Gal}({\sf N}_{\bullet}[K]/K)\] the
group of automorphisms of ${\sf N}_{\bullet}[K]$ fixing $K$;
if $L/K$ is Galois, denote by 
\[ {\sf Gal}({\sf N}_{\bullet}[L]/{\sf N}_{\bullet}[K])\] the group
of automorphisms of ${\sf N}_{\bullet}[L]$
fixing ${\sf N}_{\bullet}[K]$.

\begin{stheo2}  For all $K$, ${\sf Gal}({\sf N}_{\bullet}[K]/K)$ is trivial.
If $L/K$ is Galois, then 
\[ {\sf Gal}({\sf N}_{\bullet}[L]/{\sf N}_{\bullet}[K])
\;\;\cong\;\;
{\sf Gal}(L/K).\]
\end{stheo2}

We consider finally the case $K=\Q$.  In order to interpret the idele class group $C_{\Q}$ as a Galois group within this framework,
the operations $\oplus$ and $\otimes$ must be decoupled. 
Let us consider $\bar{\sf N}_{\bullet}[\Q ] =\PR{\sf H}_{\bullet}[\Q ]$, an abstract nonlinear number field
containing ${\sf N}_{\bullet}[\Q ]$.
 We show that
there exist flows \[ [\Phi] :\Q_{\infty}\;\longrightarrow\; 
{\sf Gal}_{\oplus}(\bar{\sf N}_{\bullet}[\Q ]/\Q) ,\quad\quad\quad\quad
[\Psi] :\Q_{\infty}\;\longrightarrow\; 
{\sf Gal}_{\otimes}(\bar{\sf N}_{\bullet}[\Q ]/\Q ) ,
\] 
where ${\sf Gal}_{\oplus}$, ${\sf Gal}_{\otimes}$
denote the groups of automorphisms of $\bar{\sf N}_{\bullet}[\Q ]$ preserving $\oplus$, $\otimes$ only.
Using the above theorem
with $L=\Q^{\rm ab}$ 
leads to the following

\begin{stheo2} There are monomorphisms
\[ {\sf C}_{\Q}\;\hookrightarrow\; 
{\sf Gal}_{\oplus}\big(\bar{\sf N}_{\bullet}[\Q^{\rm ab}]/\Q\big) ,
\quad\quad\quad\quad
{\sf C}_{\Q }\;\hookrightarrow\; 
{\sf Gal}_{\otimes}\big(\bar{\sf N}_{\bullet}[\Q^{\rm ab}]/\Q\big) .
\]
\end{stheo2}
These results are proved in \S \ref{GalGroupChap}.

\noindent {\it Remark}:  We have recently written some notes \cite{GeVeN} which expand on
the ideas described in this paper.

\vspace{5mm}

\noindent {\bf {\em Acknowledgements:}}  We thank Michael McQuillen for having pointed out several errors in the original, published version.
We would also like to thank
the International Centre for Theoretical Physics, for providing 
generous support and congenial working conditions during the period
in which this paper was written. 
This work was partially supported by {\it proyecto} 
{\bf CONACyT} {\it Sistemas Din\'amicos} G36357-E.

\section{Solenoids}\label{SolenoidChap}

We review here a notion fundamental to 
this paper.  References: \cite{Ca-Co}, \cite{Gh}, \cite{Go}, \cite{Mo-Sc}.

Let ${\sf T}$ be a $2^{\rm nd}$ countable, Hausdorff space.  An {\em n-lamination} is a $2^{\rm nd}$ countable, Hausdorff space
${\sf L}$ equipped with a maximal atlas of homeomorphisms
\[ \{ \phi_{\alpha}:{\sf L}\supset U_{\alpha}\longrightarrow
V_{\alpha}\subset\R^{n}\times {\sf T} \},\] 
in which every overlap
$\phi_{\alpha\beta}=\phi_{\beta}\circ\phi^{-1}_{\alpha}$ is of the
form
\[ \phi_{\alpha\beta}({\bf x}, {\sf t})=
( h_{\alpha\beta}({\bf x}, {\sf t}),\, f_{\alpha\beta}({\sf
t})),\] where ${\sf t}\mapsto h_{\alpha\beta}(\cdot , {\sf t})$ is a continuous family
of homeomorphisms and $f_{\alpha\beta}$ is a
homeomorphism.
We call ${\sf L}$ a {\em foliation} if ${\sf T}=\R^{k}$, a 
{\em solenoid} if ${\sf T}$ is a Cantor set.

Let $\phi$ be a chart, $D\subset\R^{n}$ an open disk,
${\sf T}'\subset {\sf T}$ open.  A {\em flowbox} is a subset of the form
$\phi^{-1}(D\times {\sf T}')$, a {\em
flowbox transversal} a subset of the form $\phi^{-1}({\bf x}\times {\sf
T}')$ and a {\em plaque} 
a subset of the form $\phi^{-1}(D\times \{ {\sf t}\})$.  A {\em leaf} is a maximal
continuation of overlapping plaques:  by definition, an
$n$-manifold.  A riemannian metric on a smooth lamination
is a family $\gamma =\{\gamma_{\ell}\}$ of smooth riemannian
metrics, one on each leaf $\ell$, which when restricted to a flowbox
gives a continuous family ${\sf t}\mapsto \gamma_{\sf t}$
of smooth metrics.
If ${\sf L}$ has the structure of a topological group such
that the multiplication and inversion maps take leaves
to leaves, ${\sf L}$ is called a {\em topological
group lamination}.  If in addition ${\sf L}$ is smooth 
and the multiplication and inversion maps are smooth along the leaves,
${\sf L}$ is called a {\em Lie
group lamination}.

Let $B$ be a manifold, $F$ a $2^{\rm nd}$ countable
Hausdorff space, $\rho : \pi_{1}B\rightarrow {\sf Homeo}(F)$ a
representation.  The quotient
\[ {\sf L}_{\rho}\;\;=\;\; \Big(\widetilde{B}\times F   \Big)\Big/ \pi_{1}B  \]
by the action $\alpha\cdot(\tilde{x}, t)=(\alpha\cdot \tilde{x},\,
\rho_{\alpha}(t))$ is a lamination called the {\em suspension} of
$\rho$.  The projection
${\sf L}_{\rho}\rightarrow B$ displays ${\sf L}_{\rho}$
as a fiber bundle with fiber
$F$, and the restriction of the projection to any leaf 
is a covering map.  

The solenoids considered in this
paper are modeled on the following example.  Let $B=\T^{d}$ = the $d$-torus and
let $F=\hat{\Z}^{d}$ = the
profinite completion
of $\Z^{d}$.  Let
\[ \rho :\pi_{1}\T^{d}\cong\Z^{d}\longrightarrow {\sf Homeo}(\hat{\Z}^{d}),\quad
\rho_{\bf
n}(\hat{\bf m})=\hat{\bf m}-{\bf n}\]
for ${\bf n}\in \Z^{d}$ and $\hat{\bf m}\in \hat{\Z}^{d}$.  The associated
suspension is called the $d$-{\em dimensional torus solenoid}
$\hat{\T}^{d}$.  Its leaves are its path
components, are dense, and may be identified with $\R^{d}$.
Moreover, $\hat{\T}^{d}$ is a compact, abelian, Lie group 
solenoid, and the additive group $\R^{d}$ sits
inside as the path-component/leaf containing the identity.

One can extend the definition of a lamination to include infinite
dimensional leaves modeled on a locally convex topological vector space $V$.  If
$V$ is a Hilbert space, we say we have a Hilbert space lamination.  
If $V=\R^{\omega}$ with the Tychonoff topology, we call it an
$\R^{\omega}$-lamination.  For Hilbert space laminations,
one may use the Fr\'{e}chet derivative to define smoothness,
for $\R^{\omega}$-laminations one uses the G\^{a}teaux
derivative.
One may make sense of all  
other notions discussed above in the infinite-dimensional setting.

\section{Adele Class Groups I: Finite Field Extensions}\label{AdClassIChap}

The material here is classical and can be found
in standard texts on algebraic number theory: \cite{Ca-Fr},
\cite{Ko}, \cite{Ne}, \cite{Ra-Va}.  We review it in order to fix notation.

Let $K$ be an algebraic number field of degree $d$ over $\Q$, 
$O_{K}$ its ring of integers.
If $K$ is Galois over $\Q$
we denote by ${\sf Gal}(K/\Q )$
the Galois group.  By a local field, we mean
a locally compact field.  
A local field of characteristic 0 is either
$\R$, $\C$ or a finite extension of $\Q_{p}$ = the field of
$p$-adic numbers.  

Let $\nu :K\rightarrow K_{\nu}$ be an embedding
where $K_{\nu}$ is a local field (necessarily of characteristic 0) 
and $\nu (K)$ is dense.  Embeddings that are related
by a continuous isomorphism of target fields different from complex conjugation
are deemed equivalent, and
an equivalence class of embedding is called a place.
When $K/\Q$ is Galois, we view
$\sigma\in {\sf Gal}(K/\Q )$ as acting on the left of the places via
$\sigma\nu:=\nu\circ\sigma$.  In practice, we shall use the word place
to mean a representative of an embedding equivalence class, 
and we will write $q_{\nu}$ for
$\nu (q)$, $q\in K$.  

When $K_{\nu}$ is isomorphic to $\R$ or $\C$, $\nu$ is said to be
a real or complex infinite place.    
If $K$ has a complex place $\mu$, then
there is an element $\sigma_{\sf conj}\in {\sf Aut}(K)$, called
$\C$-conjugation, such
that $\sigma_{\sf conj}\mu =\overline{\mu}$:
so the complex places come in conjugate pairs $(\mu ,\bar{\mu})$.
We have $d=r+2s$ where $r=$ the number of real places, $2s=$
the number of complex places. 
Let $\mathcal{P}_{\infty}=\{ \nu_{1},\dots \nu_{r},\mu_{1},\bar{\mu}_{1}, \dots ,\mu_{s},\bar{\mu}_{s}\}$,
where $\nu_{j}$, $j=1,\dots ,r$ are the real places and the
$\mu_{k},\bar{\mu}_{k}$, $k=1,\dots ,s$, are the complex places.
When $K/\Q$ is Galois, the places are either all real or
all complex, in which case we will refer to $K$ as being either real or complex.

If $\nu$ is
not infinite, it is said to be finite.
The set of finite places will be denoted $\mathcal{P}_{\sf fin}$. If
$\nu\in\mathcal{P}_{\sf fin}$, then $K_{\nu}$ has a
maximal open subring $O_{\nu}$: its ring of integers. We note that $K_{\nu}$
is locally Cantor (totally-disconnected, perfect and locally compact) and
$O_{\nu}$ is Cantor.  
 
Denote \[ K_{\infty} \;\;=\;\; \big\{  (z_{\nu})\in\C^{\mathcal{P}_{\infty}} =\; \C^{d} \; |\;\;
\bar{z}_{\nu} = z_{\bar{\nu}}
 \big\}
 \;\; \cong\;\; \R^{r}\times\C^{s}\;\;\cong\;\;\R^{d}.\] 
Note that $K_{\infty}$ has a canonical
inner-product induced from the hermitian inner product on $\C^{d}$, which decomposes as the
usual inner product
on $\R^{r}$ and the usual hermitian inner-product on $\C^{s}$.  $K$ embeds diagonally into
$K_{\infty}$ via $q\mapsto (q_{\nu})$ and the image of $O_{K}$ is
a lattice in $K_{\infty}$.   We shall identify $K$ and $O_{K}$ with their images in $K_{\infty}$.
The quotient
\[\T_{K}\;\;=\;\;K_{\infty}/O_{K}\]
is called the {\it Minkowski torus},  a torus of (real) dimension $d$.  
When $K/\Q$ is totally complex, $\T_{K}$ also has
the structure of a complex torus.

The ring of finite adeles is the restricted product 
\[ \A^{\sf fin}_{K}\;\; =\;\; 
\sideset{}{^\prime}\prod_{\nu\in\mathcal{P}_{\sf fin}}K_{\nu}\]
with respect to the $O_{\nu}$, $\nu\in\mathcal{P}_{\sf fin}$.  
By definition, this is the set of all
tuples $(q_{\nu})$ in which $q_{\nu}\in O_{\nu}$ for almost every
$\nu\in\mathcal{P}_{\sf fin}$.  $\A^{\sf fin}_{K}$ is a locally Cantor topological
ring. The adele ring is defined
\[  \A_{K} \;\;=\;\;K_{\infty}\times \A^{\sf fin}_{K}.\]
$\A_{K}$ is a locally compact ring, and a solenoid as well since it is locally
homeomorphic to (euclidean) $\times$ (Cantor).  $K$ embeds 
diagonally in $\A_{K}$ as a discrete
co-compact subgroup with respect to addition.  The quotient
\[\SO_{K}=\A_{K}/(K, +)\] is called the {\em adele class group} associated to $K$.  

Given $\mathfrak{a}$ an ideal in
$O_{K}$, denote by $\T_{\mathfrak{a}}$ the quotient
$K_{\infty}/\mathfrak{a}$, a
$d$-dimensional torus covering $\T_{K}$. 

\begin{prop}\label{invlimitwist} $\SO_{K}$ is a $d$-dimensional 
euclidean Lie group solenoid, isomorphic
to 
\begin{enumerate}
\item The inverse limit of euclidean Lie groups
\[\lim_{\longleftarrow}\;\; \T_{\mathfrak{a}},\]
where $\mathfrak{a}$ ranges over all ideals in
$O_{K}$.
\item The suspension 
\[ \big( K_{\infty}\;\times\;
\hat{O}_{K}\big)\big/ {O_{K}},\]
where $\hat{O}_{K}=\underset{\longleftarrow}{\lim}\;
O_{K}/\mathfrak{a}$.
\end{enumerate}
\end{prop}

\begin{proof}  For item (1), see \cite{Ko}, pg. 67.  Given $({\bf z},\hat{\gamma})\in K_{\infty}\times\hat{O}_{K}$, note that ${\bf z}$
defines an element $\hat{\bf z}=({\bf z}_{\mathfrak{a}})$ of 
$\underset{\longleftarrow}{\lim}\;\T_{\mathfrak{a}}$ by projection
on to each of the factors. Moreover, $\hat{\gamma}$ is by
definition a coherent sequence $\{\gamma_{\mathfrak{a}}\}$ of
deck transformations of the coverings
$\T_{\mathfrak{a}}\rightarrow\T_{K}$.  Then the association $({\bf
z},\hat{\gamma})\mapsto \hat{\gamma}\cdot\hat{\bf
z}=\big(\gamma_{\mathfrak{a}}({\bf z}_{\mathfrak{a}})\big)$ defines a
homomorphism $K_{\infty}\times\hat{O}_{K}\rightarrow \SO_{K}$
which identifies precisely the $O_{K}$-related points, and
descends to an isomorphism of $ \big( K_{\infty}\times
\hat{O}_{K}\big)\big/ {O_{K}}$ with
$\underset{\longleftarrow}{\lim}\;\T_{\mathfrak{a}}$.  
\end{proof}

By Proposition~\ref{invlimitwist}, it follows
that the path-component of 0 is a leaf canonically
isometric to $K_{\infty}$; the restriction of the projection
$\SO_{K}\rightarrow \T_{K}$ to $K_{\infty}$ is the quotient by $O_{K}$.
We will endow $\SO_{K}$ with the
riemannian metric $\rho$ along the leaves induced from the inner-product on
$K_{\infty}$. 

For example, if $K=\Q$, then $\T_{\Q}\cong S^{1}$ and $\SO_{\Q}$ is the
classical 1-dimensional solenoid obtained as the inverse limit of
circles $\R/m\Z$ under the natural covering homomorphisms.

We now suppose that $K/\Q$ is Galois and describe the actions of the
Galois groups.    If ${\bf z}=(z_{\nu})\in K_{\infty}$ and $\sigma\in {\sf Gal}(K/\Q )$
then $\sigma ({\bf z})= (z_{\sigma\nu})$.  Alternatively,                                
we can write  \[ K_{\infty}\;\;\cong\;\;
\R\otimes_{\Q}K\quad\quad\quad\quad\A_{K}\;\;\cong\;\;
\A_{\Q}\otimes_{\Q}K\] 
and the action of $\sigma$ 
on $K_{\infty}$ and on $\A_{K}$ is via $x\otimes q \mapsto
x\otimes \sigma (q)$, where $x\in \R$ or $\A_{\Q}$
as the case may be.
In any event it is clear that ${\sf Gal}(K/\Q )$ acts orthogonally.
The action of $\sigma$ on
$\A_{K}$ may be understood as the product of the actions on $K_{\infty}$ and $\A^{\sf fin}_{K}$.

Note that the image of $O_{K}$ resp. $K$ is preserved by
$\sigma$ and we induce (leafwise) isometric isomorphisms
\[\sigma:\T_{K}\longrightarrow\T_{K}\quad\quad\quad\quad \hat{\sigma}
:\SO_{K}\longrightarrow\SO_{K},\] which are intertwined by the projection
$p:\SO_{K}\rightarrow \T_{K}$ in that
$p\circ\hat{\sigma}=\sigma\circ p$. This leads to representations
\[\rho:{\sf Gal}(K/\Q )\longrightarrow{\sf Isom}(\T_{K})\quad\quad\quad\quad
\hat{\rho} :{\sf Gal}(K/\Q )\longrightarrow {\sf Isom}(\SO_{K}) ,\]
where ${\sf Isom}(\cdot )$ means the group of isometric isomorphisms.

Let $L/K$ be a finite extension of number fields.  Any place
of $L$, finite or not, defines one on $K$ by restriction. We thus obtain  
injective inclusions of vector spaces resp. rings
\[ K_{\infty}\;\hookrightarrow \;L_{\infty},\quad\quad\quad\quad
\A_{K}^{\sf fin}\;\hookrightarrow\; \A_{L}^{\sf fin},\quad\quad\quad\quad
\A_{K}\;\hookrightarrow\; \A_{L}
\] which scale inner-products/metrics by a factor of $\deg (L/K)$.
These maps in turn
induce injective homomorphisms 
\begin{equation}\label{inclusions} \T_{K}\;\hookrightarrow\; \T_{L},\quad\quad\quad\quad
\SO_{K}\;\hookrightarrow \;\SO_{L}\end{equation}
which scale metrics by $\deg (L/K)$ and
whose images are, in case $L/K$ is Galois, fixed by the action of
${\sf Gal}(L/K )$.  

We may also define maps in the opposite direction through the trace map
${\sf Tr}_{L/K}:L\rightarrow K$, which, when $L/K$ is Galois, is given by
\[ {\sf Tr}_{L/K}(\alpha )\;\; =\;\;\sum_{\sigma\in {\sf Gal}(L/K)}\sigma (\alpha ) ,\]
and which has the property that ${\sf Tr}_{L/K}(O_{L})\subset O_{K}$.
The trace map extends to a linear map ${\sf Tr}_{L/K}:L_{\infty}\rightarrow K_{\infty}$
as follows: if ${\sf Tr}_{L/K}({\bf z})={\bf w}=(w_{\nu})$ then
\[ w_{\nu}= \sum_{\nu'} z_{\nu'}  \]
where the sum is over $\nu'$ that restrict to $\nu$.  Note that
if $\nu$ is a real and $\nu'$ is complex, then complex conjugation enters into
the above sum, so that the contribution from $z_{\nu'}$ is $2{\rm Re}(z_{\nu'})$.
We thus induce epimorphisms
\[ {\sf Tr}_{L/K}:\T_{L}\rightarrow
\T_{K},\quad\quad\quad\quad\hat{\sf T}{\sf r}_{L/K}:\SO_{L}\rightarrow \SO_{K} \] 
which are, in case $L/K$ is Galois,  ${\sf Gal}(L/K )$-equivariant: ${\sf Tr}_{L/K}\circ\sigma ={\sf Tr}_{L/K}$ and
$\hat{\sf T}{\sf r}_{L/K}\circ\hat{\sigma} =\hat{\sf T}{\sf r}_{L/K}$ for all
$\sigma\in {\sf Gal}(L/K)$.  

For the extension $K/\Q$ we will write the associated trace map ${\sf Tr}={\sf Tr}_{K/\Q}$, referring
to it as the {\it absolute
trace}.  The ideal
\[  \mathfrak{d}_{K}  = \{  \alpha\in K \; |\;\; {\sf Tr}(\alpha\cdot O_{K})\subset \Z    \}^{-1} \subset O_{K} \]
is called the {\it absolute different}.

\section{Adele Class Groups II: Infinite Field Extensions}\label{AdClassIIChap}

Let $\mathcal{K}$ be a field occurring as a direct limit
$\underset{\longrightarrow}{\lim}\; K_{\lambda}$ of fields
$K_{\lambda}/\Q$ of finite degree {\em e.g.}\ 
$\bar{\Q}$ = the algebraic closure of $\Q$
or $\Q^{\rm ab}$ = the maximal
abelian extension of $\Q$.
We may associate to $\mathcal{K}$ (non locally-compact) abelian groups 
by taking the induced direct limits of inclusions (\ref{inclusions}):
\begin{equation}\label{dltorsol} \lim_{\longrightarrow} \T_{K_{\lambda}}\, ,
\quad\quad\quad\quad
\lim_{\longrightarrow}
\SO_{K_{\lambda}}. \end{equation} 
The projections $\SO_{K_{\lambda}}\rightarrow\T_{K_{\lambda}}$ induce
a projection $\underset{\longrightarrow}{\lim}\; \SO_{K_{\lambda}}\longrightarrow
\underset{\longrightarrow}{\lim}\;  \T_{K_{\lambda}}$.
If the $K_{\lambda}$ are Galois over $\Q$, then the inverse system of Galois groups $\big\{ {\sf Gal}(K_{\lambda}/\Q )\big\}$
acts compatibly on the 
direct systems of tori and solenoids,
inducing an action of the profinite Galois group
${\sf Gal}(\mathcal{K}/\Q )  =
\underset{\longleftarrow}{\lim}\; {\sf Gal}(K_{\lambda}/\Q ) $ on each of the
spaces appearing in (\ref{dltorsol}).

Consider also
\[  \lim_{\longrightarrow} (K_{\lambda})_{\infty}\, ,
\quad\quad\quad\quad  
\lim_{\longrightarrow}\A^{\sf fin}_{K_{\lambda}}\, ,
\quad\quad\quad\quad 
\lim_{\longrightarrow}\A_{K_{\lambda}},
\]
the first a (non locally-compact) topological vector space, the last two 
(non locally-compact) topological rings.
Note that we may identify
\[ \lim_{\longrightarrow}\A_{K_{\lambda}}\;\;\cong\;\;
\lim_{\longrightarrow} (K_{\lambda})_{\infty}\times\;
\lim_{\longrightarrow}\A^{\sf fin}_{K_{\lambda}}, 
\] since the direct limit maps preserve the solenoid product structure $(\text{euclidean})
\times(\text{locally Cantor})$ of the adele spaces.
There are natural inclusions 
$O_{\mathcal{K}}\hookrightarrow \underset{\longrightarrow}{\lim}\;  (K_{\lambda})_{\infty}$
and $\mathcal{K}\hookrightarrow \underset{\longrightarrow}{\lim}\; \A_{K_{\lambda}}$, and
\[ \lim_{\longrightarrow} \T_{K_{\lambda}}
\;\;\cong\;\; \Big(\lim_{\longrightarrow} (K_{\lambda})_{\infty}\Big) \Big/ O_{\mathcal{K}} 
\quad\quad\quad\quad
\lim_{\longrightarrow}
\SO_{K_{\lambda}}\;\;\cong\;\; 
\Big(\lim_{\longrightarrow}\A_{K_{\lambda}}\Big) \Big/ \mathcal{K}.\]
It follows that $\underset{\longrightarrow}{\lim}\; \SO_{K_{\lambda}}$ is a lamination,
each leaf of which may be identified with $\underset{\longrightarrow}{\lim}\;  (K_{\lambda})_{\infty}$.
Although $\underset{\longrightarrow}{\lim}\; \A^{\sf fin}_{K_{\lambda}}$
is totally disconnected and perfect, it is not locally compact 
since the direct limit
maps are not open.  Thus $\underset{\longrightarrow}{\lim}\; \SO_{K_{\lambda}}$ is
not  a solenoid.

An inner-product on the topological vector space
$\underset{\longrightarrow}{\lim}\;  (K_{\lambda})_{\infty}$ is defined by
scaling the canonical inner-product on each summand $(K_{\lambda})_{\infty}$ by 
$(\deg (K_{\lambda}/\Q))^{-1}$ and taking the direct limit.  The action of
$O_{\mathcal{K}}$ on $\underset{\longrightarrow}{\lim}\;  (K_{\lambda})_{\infty}$ preserves
this inner-product and  
we induce a riemannian metric on
$\underset{\longrightarrow}{\lim}\;  \T_{K_{\lambda}}$.
Similarly, the action of $\mathcal{K}$ on $\underset{\longrightarrow}{\lim}\; \A_{K_{\lambda}}$
is isometric along the factor $\underset{\longrightarrow}{\lim}\; (K_{\lambda})_{\infty}$, and
we induce a leaf-wise riemannian metric on 
 $\underset{\longrightarrow}{\lim}\; 
\SO_{K_{\lambda}}$.

The completion of $\underset{\longrightarrow}{\lim}\; (K_{\lambda})_{\infty}$ is
a Hilbert space denoted $\mathcal{K}_{\infty}$.  Write
\[ \A_{\mathcal{K}}\;\;=\;\;
\mathcal{K}_{\infty}\times \;\lim_{\longrightarrow}\A^{\sf fin}_{K_{\lambda}}.\]  
We define
the {\em hilbertian torus} resp.\ the {\em hilbertian adele
class group} of
$\mathcal{K}$ by
\[ \T_{\mathcal{K}}\;\;\cong\;\;  \mathcal{K}_{\infty}/ O_{\mathcal{K}}
\quad\quad\quad\quad
\SO_{\mathcal{K}}\;\cong\; \A_{\mathcal{K}}/\mathcal{K}. \]
 Thus $\SO_{\mathcal{K}}$ is a Hilbert space lamination whose leaves are isometric to
$\mathcal{K}_{\infty}$.  When the system is Galois, the Galois actions on the summands
of the direct limits preserve the scaled inner-products, and we obtain representations
\[ {\sf Gal}(\mathcal{K}/\Q )\;\longrightarrow {\sf Isom}(\T_{\mathcal{K}} )  
\quad\quad\quad\quad
{\sf Gal}(\mathcal{K}/\Q )\;\longrightarrow  {\sf Isom}(\SO_{\mathcal{K}}).\]

Any ideal $\mathfrak{a}\subset O_{\mathcal{K}}$ can be realized as the direct
limit of ideals 
$\mathfrak{a}_{\lambda}=\mathfrak{a}\cap O_{K_{\lambda}}$. 
The quotient 
$\T_{\mathfrak{a}}=\mathcal{K}_{\infty}/\mathfrak{a}$ will
be referred to as the hilbertian torus of $\mathfrak{a}$.
Consider the inverse limit
\[ \hat{O}_{\mathcal{K}}\;\; =\;\;  
\lim_{\longleftarrow} O_{\mathcal{K}}/\mathfrak{a} \]
as $\mathfrak{a}$ ranges over ideals in $O_{\mathcal{K}}$.
Since ideals in $O_{\mathcal{K}}$ 
need not have finite index,
$\hat{O}_{\mathcal{K}}$ is not compact, or even locally compact.
Nevertheless, 

\begin{prop}\label{limint}  $\hat{O}_{\mathcal{K}}\;\cong\; 
\underset{\longrightarrow}{\lim}\;  \hat{O}_{K_{\lambda}}$.
\end{prop}

\begin{proof}  Observe first that for any ideal $\mathfrak{a}\subset
O_{\mathcal{K}}$, we have $O_{\mathcal{K}}/ \mathfrak{a}\cong 
\underset{\longrightarrow}{\lim}\;  O_{K_{\lambda}}/\mathfrak{a}_{\lambda}$.
If $\mathfrak{b}\subset\mathfrak{a}$, we have a commutative diagram
\[ \begin{array}{ccccc}
  O_{\mathcal{K}}/ \mathfrak{a} & & \longleftarrow & & O_{\mathcal{K}}/ \mathfrak{b}\\
& \\
\cup & & & &\cup \\
& \\
O_{K_{\lambda}}/ \mathfrak{a}_{\lambda} & & \longleftarrow 
& & O_{K_{\lambda}}/ \mathfrak{b}_{\lambda}
\end{array}.
    \]
This allows us to interchange direct and inverse limits and write:
\[  \hat{O}_{\mathcal{K}}\;\; =\;\;
\lim_{\longleftarrow} O_{\mathcal{K}}/\mathfrak{a}\;\; \cong\;\;
\lim_{\longleftarrow}
\Big(\lim_{\longrightarrow}O_{K_{\lambda}}/\mathfrak{a}_{\lambda} \Big)\;\;\cong\;\;
\lim_{\longrightarrow} 
\hat{O}_{K_{\lambda}} \]
and the result follows.
\end{proof}

\begin{prop}\label{invlimitwistinf}  $\SO_{\mathcal{K}}$ is isomorphic to
\begin{enumerate}
\item  The inverse limit of hilbertian tori
\[ \lim_{\longleftarrow} \T_{\mathfrak{a}}\]
as $\mathfrak{a}$ ranges over ideals in $O_{\mathcal{K}}$.
\item  The suspension
\[  \Big(\mathcal{K}_{\infty}\times\hat{O}_{\mathcal{K}}\Big) 
\Big/ O_{\mathcal{K}}. \]
\end{enumerate}
\end{prop}

\begin{proof}  The spaces appearing in items (1) and (2) are isomorphic
by an argument identical to that appearing in Proposition~\ref{invlimitwist}.
That $\SO_{\mathcal{K}}$ is isomorphic to the 
inverse limit appearing in (1) follows from the same limit interchange
argument used in the proof of Proposition~\ref{limint}.
\end{proof}

Unfortunately, neither $\T_{\mathcal{K}}$ nor $\SO_{\mathcal{K}}$ are 
locally compact.
This creates complications,  since
harmonic analysis plays a fundamental role in the linking of the arithmetic
of $\mathcal{K}$ and the algebra of the Hilbert space of
$L^{2}$ functions on $\SO_{\mathcal{K}}$.
For this reason, we  consider in the next section compactifications 
coming from {\em inverse limits} of tori and solenoids.

\section{Proto Adele Class Groups}\label{ProtoAdClassChap}

Let $\mathcal{K}=\underset{\longrightarrow}{\lim}\; K_{\lambda}$ be a direct limit
of finite extensions over $\Q$.  The
trace maps induce an inverse limit
\[\hat{\mathcal{K}}\;\; =\;\;\lim_{\longleftarrow}\; K_{\lambda},\]
an abelian group with respect to addition.  This system restricts
to one of integers, however:

\begin{theo}\label{nointegers} Let $\mathcal{K}$ be a field
containing $\Q^{\rm ab}$.  Then
\[ \lim_{\longleftarrow}O_{K_{\lambda}}\; =\; \{ 0\}.\]
\end{theo}

\begin{proof} Let $\omega$ be a primitive $m$th root of unity and consider the cyclotomic
extension
$K=\Q (\omega )$.  Since $K$ is abelian, by assumption it occurs in the direct
system defining $\mathcal{K}$.  The ring of integers 
$O_{K}$ is generated by $1$ and $\omega^{j}$, where $1\leq j\leq d-1$
and $d$ is the degree of $K/\Q$. If we take say $m=2^{k}$, then $d=2^{k-1}$ and
\[ {\sf Tr}(\omega^{{j}})\;\;=\;\;\sum_{\sigma\in {\sf Gal}(K/\Q)}\sigma (\omega^{j})
\;\; =\;\; 0\] for each $j\geq 1$.  On the other hand,
${\sf Tr}(1)=d$, so  
it follows that ${\sf Tr}(O_{K})=(d )\subset\Z$.
Since $d$ can be taken arbitrarily large, this means that the only
coherent sequence that we may form in the inverse limit of rings of 
integers is $(0,0,\dots )$.  
\end{proof}

\begin{note} It is worth pointing out that normalizing
the trace map by dividing by the degree would not produce
a non-trival limit.  Indeed, if we
let $\omega$ be a primitive $p$th root of unity, $p$ a prime $>2$, then 
${\sf Tr}(O_{K})=\Z$.  Normalizing would take us out of the integers.
Thus, all that survives in the trace inverse limit is a kind of ``scaled'' 
additive number theory.
\end{note}

The trace inverse limits
\[ \hat{\T}_{\mathcal{K}}\;\; =\;\; \lim_{\longleftarrow}\;\T_{K_{\lambda}},
\quad\quad\quad\quad
\hat{\SO}_{\mathcal{K}}\;\;=\;\;\lim_{\longleftarrow}\;\SO_{K_{\lambda}}
\]
are called, respectively, the {\em proto-torus} and the
{\em proto-adele class group} of $\mathcal{K}$.  
Each is a
compact abelian group, being inverse limits of the same.  
The trace maps are natural with respect to the
epimorphisms $\SO_{K_{\lambda}}\rightarrow \T_{K_{\lambda}}$ and
induce in turn an epimorphism $\hat{\SO}_{\mathcal{K}}\rightarrow
\hat{\T}_{\mathcal{K}}$. 

Consider as well the inverse limits
\[ \hat{\mathcal{K}}_{\;\infty}\;\; =\;\; 
\lim_{\longleftarrow}\; (K_{\lambda})_{\infty}\, ,\quad\quad\quad\quad 
\hat{\A}^{\sf fin}_{\mathcal{K}}\;\; =\;\;
\lim_{\longleftarrow}\; \A^{\sf fin}_{K_{\lambda}}\, ,
\quad\quad\quad\quad 
\hat{\A}_{\mathcal{K}}\;\; =\;\;\lim_{\longleftarrow}\; \A_{K_{\lambda}}.
\] 
Observe that $ \hat{\mathcal{K}}_{\;\infty}$ is a locally convex,
(non locally compact) topological vector space whose topology
is induced from an embedding in $\R^{\omega}$. Moreover, 
$\hat{\A}^{\sf fin}_{\mathcal{K}}$ and $\hat{\A}_{\mathcal{K}}$ are
abelian topological groups only.  Note that 
\[ \hat{\A}_{\mathcal{K}}\;\;\cong\;\; \hat{\mathcal{K}}_{\;\infty}
\times \hat{\A}^{\sf fin}_{\mathcal{K}}. \]
The space $\hat{\A}^{\sf fin}_{\mathcal{K}}$ is totally-disconected,
perfect but not locally compact.  
There is a natural inclusion 
$\hat{\mathcal{K}}\;\hookrightarrow \;\hat{\A}_{\mathcal{K}}$
and
\[ \hat{\SO}_{\mathcal{K}}\;\;\cong\;\;
\hat{\A}_{\mathcal{K}}/\hat{\mathcal{K}}.
\]
Since the action of $\hat{\mathcal{K}}$
locally preserves the product structure, $\hat{\SO}_{\mathcal{K}}$ is an
$\R^{\omega}$-lamination.
When the system is Galois, the trace system maps  
are compatible with the Galois
actions and we obtain representations of ${\sf Gal}(\mathcal{K}/\Q )$
on $\hat{\T}_{\mathcal{K}}$ and on 
$\hat{\SO}_{\mathcal{K}}$, acting by (leaf-preserving) topological isomorphisms.

\begin{theo}  Let $\mathcal{K}=\underset{\longrightarrow}{\lim}\;  K_{\lambda}$
be an infinite field extension
containing $\Q^{\rm ab}$.  Then the covering homomorphisms 
$(K_{\lambda})_{\infty}\rightarrow \T_{K_{\lambda}}$ induce a continuous
{\rm monomorphism}
\[ \hat{\mathcal{K}}_{\;\infty}\;\;\longrightarrow\;\;
\hat{\T}_{\mathcal{K}}.\] 
\end{theo}

\begin{proof}  The inverse limits giving rise to each of 
$\hat{\mathcal{K}}_{\;\infty}$ and
$\hat{\T}_{\mathcal{K}}$
also give rise to a continuous homomorphism 
$\hat{\mathcal{K}}_{\;\infty}\rightarrow
\hat{\T}_{\mathcal{K}}$.  An element in the kernel 
must be a coherent sequence of algebraic integers {\em i.e.}\ an element
of $\underset{\longleftarrow}{\lim}\;  O_{K_{\lambda}}$.  By Theorem~\ref{nointegers}, the latter 
is trivial and the map is an injective.  \end{proof}

The map $\hat{\mathcal{K}}_{\;\infty}\rightarrow
\hat{\T}_{\mathcal{K}}$ is not surjective since the fibers of
the projections $(K_{\lambda})_{\infty}\rightarrow \T_{K_{\lambda}}$
are not compact.  Nonetheless, the image of $\hat{\mathcal{K}}_{\;\infty}$ is dense
in the compact  $\hat{\T}_{\mathcal{K}}$.  Thus $\hat{\T}_{\mathcal{K}}$ is laminated by
the cosets of the image of $\hat{\mathcal{K}}_{\;\infty}$.  In fact,

\begin{theo}\label{torisosol}  Let $\mathcal{K}=\underset{\longrightarrow}{\lim}\;  K_{\lambda}$
be an infinite field extension
containing $\Q^{\rm ab}$.  Then the trace map inverse systems induce an
{\rm isomorphism}
\[ \hat{\SO}_{\mathcal{K}}\;\;\longrightarrow\;\;
\hat{\T}_{\mathcal{K}}\]
of topological groups. 
\end{theo}

\begin{proof}  For each $K_{\lambda}$, the kernel of 
$\SO_{K_{\lambda}}\rightarrow \T_{K_{\lambda}}$ is
$\hat{O}_{K_{\lambda}}$, hence the kernel of $\hat{\SO}_{\mathcal{K}}\rightarrow
\hat{\T}_{\mathcal{K}}$ is $\underset{\longleftarrow}{\lim}\; \hat{O}_{K_{\lambda}}=0$.
Since the fibers of the projections $\SO_{K_{\lambda}}\rightarrow \T_{K_{\lambda}}$
are compact, the induced map of limits $\hat{\SO}_{\mathcal{K}}\rightarrow
\hat{\T}_{\mathcal{K}}$ is surjective.
Since each map $\SO_{K_{\lambda}}\rightarrow \T_{K_{\lambda}}$ is open with compact fibers, the limit map
$\hat{\SO}_{\mathcal{K}}\rightarrow
\hat{\T}_{\mathcal{K}}$ is also open.
\end{proof}

When $\Q^{\rm ab}\subset\mathcal{K}$,
there are no exact analogues of Propositions~\ref{invlimitwist}
and \ref{invlimitwistinf}, due to the lack of a notion
of integers in $\hat{\mathcal{K}}$.  On the other hand,
given any 
subfield $K\subset\mathcal{K}$ 
of finite degree over $\Q$, one can find ``level-$K$''
suspension structures on $\hat{\T}_{\mathcal{K}}\cong\hat{\SO}_{\mathcal{K}}$.  Specifically, it
is not difficult to see that the kernel of the projection 
$\hat{\T}_{\mathcal{K}}\rightarrow\T_{K}$ is isomorphic to
\[  \hat{O}_{K}\times \hat{\T}_{\mathcal{K}/K},\]
where $\hat{\T}_{\mathcal{K}/K}$ is the inverse limit of tori
occurring as the connected components of $0$ of the kernels of the projections
$\T_{L}\rightarrow \T_{K}$.  Then
\[ \hat{\T}_{\mathcal{K}}
\;\;\cong\;\; 
\Bigg( K_{\infty}\times
\Big( \hat{O}_{K}\times \hat{\T}_{\mathcal{K}/K}  \Big) \Bigg)\Bigg/ O_{K}.\]
These representations are related by homeomorphisms induced by 
${\sf Tr}: L\rightarrow K$, but they do not survive the trace inverse
limit.  Similar ``level-$K$'' suspension representations are available
for $\hat{\SO}_{\mathcal{K}}$ as well.

The relationship between the proto constructions of
this chapter with the hilbertian constructions of the previous chapter
is described by the following

\begin{theo}  There are canonical 
inclusions
\[ \T_{\mathcal{K}}\;\hookrightarrow\;\hat{\T}_{\mathcal{K}},\quad\quad\quad\quad
\mathcal{K}_{\;\infty}\;\hookrightarrow\hat{\mathcal{K}}_{\;\infty},\quad\quad\quad\quad 
\SO_{\mathcal{K}}\;\hookrightarrow\;  \hat{\SO}_{\mathcal{K}}\]
with dense images, which are ${\sf Gal}(\mathcal{K}/\Q )$-equivariant in case
the defining system is Galois.
\end{theo}

\begin{proof}  We begin by
defining the inclusion 
$\underset{\longrightarrow}{\lim}\;  (K_{\lambda})_{\infty}\;\hookrightarrow\hat{\mathcal{K}}_{\;\infty}$.
Let ${\bf z}_{\lambda_{0}}\in 
(K_{\lambda_{0}})_{\infty}$ and define
${\bf z}_{\lambda_{0}}\mapsto ({\bf z}_{\lambda})$ as follows.  If $\lambda$ is an index below 
$\lambda_{0}$,
project ${\bf z}_{\lambda_{0}}$ by the appropriate trace map to get the ${\bf z}_{\lambda}$ coordinate.
If $\lambda$ is above $\lambda_{0}$, we map ${\bf z}_{\lambda_{0}}$ into $(K_{\lambda})_{\infty}$ by the inclusion
$  (K_{{\lambda}_{0}})_{\infty}\hookrightarrow (K_{\lambda})_{\infty}$ and scale by $1/d$ where $d$ is the degree
of $K_{\lambda}/K_{\lambda_{0}}$.  This prescription defines a coherent sequence hence an element of 
$\hat{\mathcal{K}}$.  
This map clearly defines an injective homomorphism of vector spaces
$\underset{\longrightarrow}{\lim}\; (K_{\lambda})_{\infty}\hookrightarrow\hat{\mathcal{K}}$.  
By virtue
of the scaling, this map is also continuous with regard to the inner-product
defining $\mathcal{K}_{\;\infty}$.  Now let $\{ {\bf z}_{\lambda_{i}}\}
\subset\underset{\longrightarrow}{\lim}\; (K_{\lambda})_{\infty}$ be a Cauchy sequence defining
an element of $\mathcal{K}_{\;\infty}$.  For every index $\beta$, if
we let ${\sf Tr}_{\lambda_{i},\beta }$ denote the trace map from
$(K_{\lambda_{i}})_{\infty}$ to $(K_{\beta})_{\infty}$, then
$\big\{{\sf Tr}_{\lambda_{i},\beta }({\bf z}_{\lambda_{i}})\big\}$ is Cauchy in 
$(K_{\beta})_{\infty}$.  But this means precisely that 
the image of $\{ {\bf z}_{\lambda_{i}}\}$ in $\hat{\mathcal{K}}_{\;\infty}$
is convergent.
Thus the inclusion 
$\underset{\longrightarrow}{\lim}\;  (K_{\lambda})_{\infty}\;
\hookrightarrow\hat{\mathcal{K}}_{\;\infty}$
extends to $\mathcal{K}_{\;\infty}$.
The other inclusions
are induced by this one.  That the images are dense follows from the definition of the above map:
given any element of 
$\hat{\mathcal{K}}_{\;\infty}$ and any level $\lambda$, 
there is an element of $\mathcal{K}_{\;\infty}$
agreeing up to level $\lambda$.
\end{proof}

\section{Hyperbolizations}\label{hyperbolizations}

Let $K/\Q$ be finite.  For each real place $\nu$ one pairs $K_{\nu}\cong\R$ with a single factor of
$(0,\infty )$ to form the upper half-plane factor $\HP_{\nu}=\R\times i(0,\infty )$, which
we equip with the hyperbolic metric in the usual way.  

For a complex place pair $(\mu ,\bar{\mu})$, we must alter slightly this prescription
in order to take into account the complex algebra of the factor 
\[  \C_{(\mu ,\bar{\mu})} = \{  (z,\bar{z}) \; |\;\; z\in\C \}\cong \C . \]
In order to do this, it is convenient to regard the factor $(0,\infty )\times (0,\infty )$
to be attached to  $\C_{(\mu ,\bar{\mu})}$
 as being {\it complex}.
Thus let $\B$ be the complex quarter space
\[ \B  = \{ b=s+it\in\C\; |\;\; 0< s,t\}\]  and define 
\[  \HP_{(\mu ,\bar{\mu})}  = \{ (z,\bar{z})\times (b,-\bar{b}) \; |\;\;
z\in \C, \; b\in  \B \}\subset (\C\times \B )\times (\bar{\C}\times (-\bar{\B} ))\subset\C^{2}\times\C^{2}. \]

Let $-i\HP^{2}$ denote the right half-plane. 
Then we may identify $\HP_{(\mu ,\bar{\mu})}$ with $\HP^{2}\oplus -i\HP^{2}$
via the map
\[    (z,\bar{z})\times (b,-\bar{b})\mapsto\frac{1}{2} ( z+\bar{z} +  b-\bar{b},  \, z-\bar{z}+ b+\bar{b})
= (x+it, s+iy) \equiv (u,v)
  \]
for $z=x+iy$.  We give $\HP_{(\mu ,\bar{\mu})}$ the product hyperbolic metric coming from this identification, and also write
\begin{equation}\label{upperright} \HP_{(\mu ,\bar{\mu})}=\HP_{\mu}\oplus -i\HP_{\mu} \end{equation}
so as to be able to refer to the upper and right half-planes associated to $(\mu ,\bar{\mu} )$.
In the sequel, a function defined in $ \HP_{(\mu ,\bar{\mu})}$ will be considered holomorphic
if it is holomorphic in each of the variables $(u,v)=(x+it, s+iy)$ separately.  

Notice that the complex conjugation acting in $\C_{(\mu ,\bar{\mu})}$  extends to 
$\HP_{(\mu ,\bar{\mu})}$ by the identity in the $b$ variable, and that in the
$(u,v)$ variables, acts via 
\[ (u,v)\mapsto (u,\bar{v}), \]
thus defining an (orientation-reversing) isometry of $\HP_{(\mu ,\bar{\mu})}$.

Finally, we define the hyperbolization 
\[ \HP_{K}\;\; =\;\; \HP_{K}^{\R}\times   \HP_{K}^{\C} \;\; := \;\;  \prod \HP_{\nu} \times \prod \HP_{(\mu ,\bar{\mu})}\;\; =\;\;
 \prod \HP_{\nu} \times \big(\prod \HP_{\mu}\oplus -i\HP_{\mu}\big)
\;\; \cong \;\;K_{\infty}\times (0,\infty)^{d} .\]
Thus $\HP_{K}$ has the structure of a $d$-dimensional 
complex polydisk equipped with
the product riemannian metric.
For $\nu$ a real place we write $\tau_{\nu} = x_{\nu}+it_{\nu}$ for a point
of $\HP_{\nu}$, and for $\mu$ complex
we write  
\[ (\kappa_{\mu}, \bar{\kappa}_{\mu}) := ((z_{\mu},b_{\mu}), (\bar{z}_{\mu},-\bar{b}_{\mu}))  \cong (u_{\mu},v_{\mu}).\]
We write points of $\HP_{K}$ in the form ${\boldsymbol \rho}={\boldsymbol \tau}\times {\boldsymbol \kappa}$
where ${\boldsymbol \tau} = (\tau_{1}, \dots, \tau_{r})\in  \HP_{K}^{\R}$ are the coordinates of the ``real hyperbolization'' and
${\boldsymbol \kappa} = ((\kappa_{1}, \bar{\kappa}_{1}), \dots , (\kappa_{s}, \bar{\kappa}_{s}))\in \HP_{K}^{\C} $ are those of the ``complex
hyperbolization".
When $K/\Q$ is Galois we have $\HP_{K}= \HP_{K}^{\R}$ or $ \HP_{K}^{\C}$
depending on whether $K$ is real or complex, and
the action of ${\sf Gal}(K/\Q )$ extends to an isometric action
on $\HP_{K}$ by acting trivially in the extended coordinates.  

The subgroups $O_{K}$ and $K$ of
$K_{\infty}$, viewed as groups of translations, extend to translations
of $\HP_{K}$ that are isometries.  The quotients
\begin{equation}\label{hyperbolized}
\mathfrak{T}_{K} \;\; =\;\;  \HP_{K}/O_{K}\;\;\cong\;\; 
(\Delta^{\ast})^{d}, \quad\quad\quad\quad
\hat{\mathfrak{S}}_{K}  \;\; =\;\;  
\Big( \HP_{K}\times \A^{\sf fin}_{K}\Big) \Big/ K \;\; \approx\;\;
\SO_{K}\times (0,\infty )^{d}
\end{equation}
are referred to as the {\em hyperbolized torus} and the {\em hyperbolized adele
class group} 
of $K$. (In the above, $\Delta^{\ast}$ denotes the punctured hyperbolic disk.) 
For $L/K$ a finite extension of finite degree extensions of $\Q$,
the canonical inclusions $\mathfrak{T}_{K}\hookrightarrow\mathfrak{T}_{L}$
and $\HSO_{K}\hookrightarrow\HSO_{L}$ are isometric up to
the scaling factor $\deg (L/K)$.  If $L/K$ is Galois,
the action of ${\sf Gal}(L/K)$ on $\T_{L}$ and $\SO_{L}$ extends 
by isometries to both of $\mathfrak{T}_{L}$ and $\HSO_{L}$ restricting
to the identity on the images of $\mathfrak{T}_{K}$
and $\HSO_{K}$.

In the case of an infinite field extension $\mathcal{K}$, 
one follows the prescription of the preceding paragraphs using the hilbertian
torus and solenoid $\T_{\mathcal{K}}$
and $\SO_{\mathcal{K}}$.  Thus $\HP_{\mathcal{K}}$
is an
infinite product of hyperbolic planes.  We obtain an infinite-dimensional 
hyperbolized torus $\mathfrak{T}_{\mathcal{K}}$ and hyperbolized solenoid
$\HSO_{\mathcal{K}}$ upon quotient by $O_{\mathcal{K}}$ and $\mathcal{K}$.

\section{The Character Field}\label{CharFieldChap}

For $G$ a locally compact abelian group, by a character we mean a
continuous homomorphism $\chi :G\rightarrow U(1)$, where
$U(1)\subset \C$ is the unit circle.  An operation $\oplus$
on characters is defined by multiplying their values, and the set of characters
\[  {\sf Char}(G)\;\;=\;\;{\sf Hom}_{\rm cont}(G,U(1))\]
is itself a locally compact abelian group.
In the case of interest, the projection $\SO_{K}\rightarrow\T_{K}$ induces
an inclusion 
\[  {\sf Char}(\T_{K})\hookrightarrow {\sf Char}(\SO_{K}) \]
and so we will view $ {\sf Char}(\T_{K})$ as a subgroup of ${\sf Char}(\SO_{K})$.
If $\mathcal{K}=\underset{\longrightarrow}{\lim}\;  K_{\lambda}$ is an infinite
extension over $\Q$, we have a similar inclusion
\[ {\sf Char}(\hat{\T}_{\mathcal{K}})\hookrightarrow {\sf Char}(\hat{\SO}_{K})\]
and corresponding convention.

The purpose of this section is to give a proof of the well-known identification
${\sf Char}(\SO_{K}) \cong (K,+)$, and show that this identification may be
used to view ${\sf Char}(\SO_{K})$ as a field, in a way which is natural with respect to the trace maps.  

\begin{lemm}\label{tracelemma}  Let $K/\Q$ be a finite extension.  If $\;{\bf w}\in K_{\infty}$
has the property that \[ {\sf Tr}({\bf w} K)\subset \Q,\] then ${\bf w}\in K$.
\end{lemm}

\begin{proof}  First suppose that $K$ is Galois over $\Q$.  Let $\alpha_{1},\dots ,\alpha_{d}$ be an integral basis of $K$, and
let $A$ be the $d\times d$ invertible matrix whose $ij$-element is $\nu_{j}(\alpha_{i})$ where
the $\nu_{1},\dots ,\nu_{d}$ are the places of $K$.  Then the hypothesis
${\rm Tr}({\bf w} K)\subset \Q$ implies that
$A{\bf w}={\bf q}\in \Q^{d}$ or 
${\bf w}\in A^{-1}\Q^{d}$.  Since $K$ is Galois, it is normal, hence all of the entries
of $A$ belong to $K$ (see \cite{La}),  thus the coordinates of ${\bf w}$ belong to $K$.  Suppose
that nevertheless ${\bf w}\not\in K$.  Then there exists some element  $\sigma_{0}\in{\sf Gal}(K/\Q )$
and a coordinate $w_{\nu_{0}}$ such that $\sigma_{0} (w_{\nu_{0}} )\not= w_{\sigma_{0}\nu_{0}}$.  Let $A^{\nu}$
denote the column indexed by a place $\nu$ ({\it i.e.}\ the vector $(\nu (\alpha_{1} ), \dots , 
\nu (\alpha_{d} ))^{T}$) so that
\begin{equation}\label{vectoreqn}
   \sum w_{\nu}A^{\nu} = {\bf q}\in \Q^{d} . 
   \end{equation}
Acting by $\sigma_{0}$ on this equation fixes ${\bf q}$, permutes the column vectors $A^{\nu}$,  $\sigma_{0} (A^{\nu}) =A^{\sigma_{0}\nu}$, 
but does not similarly permute
the entries of ${\bf w}$.  Thus there exists a vector ${\bf w}'\not= {\bf w}$, defined $w'_{\sigma_{0} \nu} =\sigma_{0} (w_{\nu})$, for which $A{\bf w}' = {\bf q}$, implying that the kernel of $A$ is nontrivial (since $A({\bf w}'-{\bf w})=0$ and  ${\bf w}'-{\bf w}\not=0$), contradiction.

Now suppose that $K/\Q$ is not Galois.  As in the previous paragraph, fix the basis $\alpha_{1},\dots ,\alpha_{d}$ and note that
 the columns of the embedding matrix
$A$ continue to satisfy the vector equation (\ref{vectoreqn}).   Let $L/\Q$ be a finite Galois extension containing the images of $K$ by its places: for example, if one writes
$K=\Q (\alpha )$ one can take $L$ to be the splitting field of the minimal polinomial of $\alpha$.
Note then that $L/K$ is Galois.
Then $A$ has entries belonging to $L$, as does its inverse, so
that each coordinate of ${\bf w}$ belongs to $L$.
Consider the usual diagonal embedding of vector spaces $\iota :K_{\infty}\hookrightarrow L_{\infty}$ defined by 
\[ \iota ({\bf w})_{\mu} = w_{\nu} \]
for each $\mu$ an $L$-place with restriction $\mu |_{K}=\nu$.  In addition, we have the scaled embedding
 $\bar{\iota} :K_{\infty}\hookrightarrow L_{\infty}$,
\[ \bar{\iota}({\bf w})_{\mu} = \frac{w_{\nu}}{{\rm mult}(\nu )},  \] where $\mu |_{K}=\nu$ and ${\rm mult}(\nu)$ is the number of places having restriction to $K$ equal to $\nu$.
The basis 
$\alpha_{1},\dots ,\alpha_{d}$ defines a $d\times (d[L:K])$ matrix, denoted $B$, whose $\mu$th-column is 
$B^{\mu}=(\bar{\iota}(\alpha_{1} )_{\mu }, \dots , 
(\bar{\iota}(\alpha_{d} )_{\mu })^{T}$ (here we are identifying the basis elements with vectors
in $K_{\infty}$).
Since $A{\bf w}={\bf q}\in \Q^{d}$, we have $B\iota ({\bf w}) ={\bf q}\in \Q^{d}$
 and we have the analogue of the vector equation (\ref{vectoreqn}):
\[
   \sum \iota({\bf w})_{\mu}B^{\mu} = {\bf q}\in \Q^{d} . 
\]
Note that if $\mu$ and $\mu'$ have the same restriction $\nu$ to $K$ then $B^{\mu}=B^{\mu'}$.

Suppose that ${\bf w}$ does not belong to $K$.   
Then $\iota ( {\bf w})$ does not belong to $L$: for if
$\iota ( {\bf w})=\beta \in L$, then by definition of $\iota$, for all
$\mu$ extending the identity place on $K$, we must have $\mu (\beta ) =\beta$.
But the places extending the identity place on $K$ comprise the Galois group
 ${\sf Gal}(L/K)$, hence
 $\beta$ belongs
to the fixed field of ${\sf Gal}(L/K)$ i.e.\ $\iota ( {\bf w})\in K$ viewed as a subfield
of $L_{\infty}$.  But then this would imply that ${\bf w} \in K$ contrary to our hypothesis.

Since $\iota ( {\bf w})$ does not belong to $L$,  there exists $\sigma_{0}\in {\sf Gal}(L/\Q )$ and a coordinate $\mu_{0}$
such that \[ \sigma_{0} (\iota({\bf w})_{\mu_{0}})\not=\iota ( {\bf w})_{\sigma_{0}\mu_{0} }.\]  Then as before
we deduce a vector ${\bf y}$, defined $y_{\sigma_{0}\mu}:= \sigma_{0} (\iota ({\bf w})_{\mu})$, and 
for which $B{\bf y}={\bf q}$.  This vector ${\bf y}$ satisfies $y_{\mu_{1}}=y_{\mu_{2}}$
whenever $\mu_{1}$ and $\mu_{2}$ have the same restriction $\nu$ to $K$,
so that it
must belong to the image of the diagonal embedding $\iota$.
  Thus there exists a vector ${\bf w}'\in K_{\infty}$ distinct from ${\bf w}$
and which satisfies $A{\bf w}' = {\bf q}$, again contradicting the invertibility of $A$.
 
\end{proof}

Recall that the inverse different of a finite degree algebraic number field $K/\Q$ is the 
$O_{K}$-module
\[   \mathfrak{d}_{K}^{-1} = \{ \alpha\in K|\; {\rm Tr}_{K/\Q}(\alpha O_{K})\subset \Z\} \supset O_{K}.\]
Note that if $L/K$ is finite and $\alpha\in  \mathfrak{d}_{K}^{-1} $ then for all $\beta\in O_{L}$,
\[ {\rm Tr}_{L/\Q}(\alpha \beta ) = {\rm Tr}_{K/\Q}(\alpha \cdot {\rm Tr}_{L/K}(\beta )) \subset\Z \]
(since ${\rm Tr}_{L/K}(\beta )\subset O_{K}$)
so that 
$   \mathfrak{d}_{K}^{-1} \subset  \mathfrak{d}_{L}^{-1}$. 
Then if $\mathcal{K}/\Q$ is an infinite degree algebraic extension, we define
\[    \mathfrak{d}_{\mathcal{K}}^{-1} = \lim_{\longrightarrow}  \mathfrak{d}_{K}^{-1}  \]
where the $K\subset\mathcal{K}$ range over finite subextensions of $\Q$.

\begin{theo}\label{Pontryagin} Suppose that 
\begin{enumerate}  
\item[\fbox{\bf a}]  $K$ is a finite field extension over $\Q$.  Then
${\sf Char}(\SO_{K})$ possesses
a second operation $\otimes$ making it a field
and for which their is a canonical isomorphism
\[ {\sf Char}(\SO_{K})\;\;\cong\;\; K   \]
which is natural with respect to the trace maps. This isomorphism
identifies ${\sf Char}(\T_{K})$ with the inverse different $\mathfrak{d}_{K}^{-1}$ and in particular,
there is a canonical embedding of the ring $O_{K}$ in ${\sf Char}(\T_{K})$.
\item[\fbox{\bf b}] $\mathcal{K}=\underset{\longrightarrow}{\lim}\;  K_{\lambda}$ is an infinite
extension over $\Q$.  Then  ${\sf Char}(\hat{\SO}_{\mathcal{K}})$ possesses
a second operation $\otimes$ making it a field
and for which their is a canonical isomorphism
\[  {\sf Char}(\hat{\SO}_{\mathcal{K}})
\;\;\cong\;\; \mathcal{K}  \]
which is natural with respect to the trace maps.  
This isomorphism
identifies ${\sf Char}(\hat{\T}_{\mathcal{\mathcal{K}}})$ with the inverse different $\mathfrak{d}_{\mathcal{\mathcal{K}}}^{-1}$.  In particular,
there is a canonical embedding of the ring $O_{\mathcal{\mathcal{K}}}$ in ${\sf Char}(\hat{\T}_{\mathcal{\mathcal{K}}})$.

\end{enumerate}
\end{theo}

\begin{note}  The statement that ${\sf Char}(\SO_{K})\cong (K,+)$ is classical and has been known since the time of Tate's thesis \cite{Ta}.  
\end{note}

\begin{proof}  We begin with {\small \fbox{\bf a}}.  Consider the character on $K_{\infty}$ defined
\[  \phi_{\infty}^{K}({\bf z})= \exp (2\pi i  \, {\sf Tr}({\bf z})). \]
Note that $ \phi_{\infty}^{K}$ extends continuously to
a character $\phi^{K}$ of
$\SO_{K}$.   Indeed, $ \phi_{\infty}^{K}$ is invariant w.r.t.\ translation by $O_{K}$ and so induces
a character $\T_{K}\rightarrow {\sf U}(1)$, which, when composed with the projection
$\SO_{K}\rightarrow \T_{K}$, defines the extension $\phi^{K}$.  Note that for any finite
degree extension $L/K$ we have $\phi^{L}= \phi^{K}\circ \hat{\sf Tr}_{L/K}$.
For each $\alpha\in K$ define $\phi_{\alpha}^{K}$ by 
\[ \phi_{\alpha}^{K}(\hat{\bf z})=
\phi^{K}(\alpha\hat{\bf z})\] (here $\alpha\hat{\bf z}$ is defined for all $\hat{\bf z}\in \SO_{K}$ since $\SO_{K}$ is a $K$-vector space).
Then the map 
\begin{equation}\label{charmap}
\alpha\mapsto \phi_{\alpha}^{K}
\end{equation}
 defines a monomorphism 
$(K,+)\rightarrow {\rm Char}(\SO_{K})$.  We will show that
it is an isomorphism.  So let $\phi :\SO_{K}\rightarrow U(1)$ be any character.  Then the restriction
to the leaf through $0$, which is just $K_{\infty}$, is of the form 
\[
{\bf z}\mapsto \exp (2\pi i {\sf Tr}({\bf w}{\bf z}))
\]
for
some ${\bf w}\in K_{\infty}$.  On the other hand, we may also restrict $\phi$ to the
transversal through $0$, which is just
$\hat{O}_{K}\subset \SO_{K}$, where it must factor through some finite quotient 
$O_{K}/\mathfrak{a}$.
It follows that $\phi$ is the pullback of a  character
\[ \T_{\mathfrak{a}}=K_{\infty}/\mathfrak{a}\longrightarrow U(1) \]
for some ideal $\mathfrak{a}\subset O_{K}$.
In particular, ${\sf Tr}({\bf w}\mathfrak{a})\subset \Z$ which implies that
${\sf Tr}({\bf w}K)\subset \Q$, so by Lemma \ref{tracelemma}, ${\bf w}=\alpha\in K$
for some $\alpha$ 
and $\phi = \phi_{\alpha}^{K}$.
Thus the map $\alpha\mapsto \phi_{\alpha}^{K}$ defines an isomorphism between
$(K,+)$ and
${\sf Char}(\SO_{K})$.  One then pushes forward
the product operation of $K$ to define $\otimes$.  
Naturality is an expression of the commutative diagram
\[ 
\begin{CD}
L    @>\cong>>    {\rm Char}(\SO_{L})\\
\cup @.                 @AA\hat{\sf T}{\sf r}_{L/K}^{\ast}A\\
K    @>\cong>>    {\rm Char}(\SO_{K})
\end{CD}
\]
whose commutativity follows from the fact that for all $\alpha\in K$,
\[ \hat{\sf T}{\sf r}_{L/K}^{\ast} (\phi^{K}_{\alpha}) = \phi^{K}_{\alpha}\circ  \hat{\sf T}{\sf r}_{L/K} =
\phi^{L}_{\alpha}.
  \]
The map (\ref{charmap}) identifies ${\sf Char}(\T_{K})$
with $\mathfrak{d}_{K}^{-1}\supset O_{K}$ so that the restriction of $\otimes$ to ${\sf Char}(\T_{K})$ makes
the latter an $O_{K}$-module isomorphic to $\mathfrak{d}_{K}^{-1}$.

To prove {\small \fbox{\bf b}} , recall that if $L/K$ is a finite extension of number fields
finite over $\Q$, then the trace projection $\hat{\sf T}{\sf r}_{L/K}:\SO_{L}\rightarrow \SO_{K}$
induces an inclusion $\hat{\sf T}{\sf r}_{L/K}^{\ast}:{\sf Char}(\SO_{K})\hookrightarrow {\sf Char}(\SO_{L})$
of fields.  Then the direct limit
\begin{equation}\label{fldlimit}  \lim_{\longrightarrow}{\sf Char}(\SO_{K_{\lambda}}) \end{equation}
is a field isomorphic to $ \mathcal{K}$.  The direct limit (\ref{fldlimit})
has dense image in ${\sf Char}(\hat{\SO}_{\mathcal{K}})$, but the latter is discrete
since $\hat{\SO}_{K}$ is compact, so they are equal.  
Thus every character $\phi$ is the pull-back of one
defined on $\SO_{K}$ where $K/\Q$ is of finite degree.
 Again, this isomorphism is trace natural. It
 identifies ${\sf Char}(\hat{\T}_{\mathcal{K}})$
with $\mathfrak{d}_{\mathcal{K}}^{-1}$ and so as in the finite extension case, the former
has a module structure over the subring corresponding to $O_{\mathcal{K}}$.
\end{proof}

Note that for $K=\Q$, we have ${\sf Char}(\T_{\Q})$ is a ring since $\mathfrak{d}_{K}=\Z$.  On the other hand
since $\mathfrak{d}_{K}\not = O_{K}$ for $K\not=\Q$, we see that ${\sf Char}(\T_{K})$ is strictly larger than $O_{K}$: it contains
$O_{K}$ as a ring, and is to be viewed as an $O_{K}$-module extension of $O_{K}$.  If $\mathcal{K}\supset \Q^{\rm ab}$, then $\mathfrak{d}_{\mathcal{K}}^{-1}=\mathcal{K}$ and 
${\sf Char}(\hat{\T}_{\mathcal{K}})={\sf Char}(\hat{\SO}_{\mathcal{K}})$, which was anticipated
by Theorem \ref{torisosol}.

 \section{Graded Holomorphic Functions}\label{GradedChap}

Let $\mu$ be the unit mass Haar measure on $\SO_{K}$ and
let $L^{2}(\SO_{K},\C )$ 
be the associated space of square integrable complex
valued functions on $\SO_{K}$.  The characters $\{\phi_{\alpha}\}$ form a complete
orthonormal system in $L^{2}(\SO_{K},\C )$ and so every element 
$f\in L^{2}(\SO_{K},\C )$
has the development
\[  f\;\;=\;\; \sum a_{\alpha}\phi_{\alpha},  \]
where $\{ a_{\alpha}\} \in l^{2}(K )$, $\phi_{\alpha}\in {\sf Char} (\SO_{K})$
and equality is taken w.r.t.\ the $L^{2}$ norm. 
We note that since we are in a Hilbert space, by Parseval's Lemma, we have
$\| f\|^{2}=\sum |a_{\alpha}|^{2} $ and so the sum
converges with respect to any well-ordering of the indexing set $K$. 
The space $L^{2}(\T_{K},\C )$ may be identified with the subspace of
$L^{2}(\SO_{K},\C )$ whose Fourier series satisfy $a_{\alpha}=0$ for $\alpha\notin \mathfrak{d}^{-1}_{K}$.
If we restrict $f$ to the dense leaf
$K_{\infty}$, then we may identify 
$\phi_{\alpha}({\bf z})=\exp (2\pi i{\rm Tr}(\alpha{\bf z}))\equiv {\boldsymbol \zeta}^{\alpha}$ and
write $f$ in the form of an $L^{2}$ Puiseux series
\[ f( {\boldsymbol \zeta})\;\; =\;\; \sum a_{\alpha} {\boldsymbol \zeta}^{\alpha}  .\]

Let $f,g\in L^{2}(\SO_{K},\C )$ be given by the developments
$f=\sum a_{\alpha}\phi_{\alpha}$, $g=\sum b_{\alpha}\phi_{\alpha}$.
The {\em Cauchy product} $f\oplus g$ is defined as the $L^{2}$
extension of the point-wise product
of continuous functions, that is,
\[  f\oplus g\;\; =\;\;\sum_{\alpha} c_{\alpha}\phi_{\alpha} \;\; =
\;\;\sum_{\alpha}\left(\sum_{\alpha_{1}+\alpha_{2}=\alpha }a_{\alpha_{1}}b_{\alpha_{2}}\right) \phi_{\alpha} ,\]
provided that the sum on the right converges.  In this regard, note that $c_{\alpha}$
is equal to the inner-product $\langle f, \sigma_{\alpha}\bar{g}\rangle$, where $\sigma_{\alpha}\bar{g}=\sum_{\beta\in K} \bar{b}_{\alpha-\beta}\phi_{\beta}\in L^{2}(\SO_{K},\C )$, hence
is an unambiguously defined complex number.  
Thus the Cauchy product is defined whenever the sequence $\{ c_{\alpha}\}$
belongs to $l^{2}(K )$.
Given 
$f\in L^{2}(\SO_{K},\C )$, we denote by ${\sf Dom}_{\oplus}(f)$ the
set of $g\in L^{2}(\SO_{K},\C )$ for which $f\oplus g$ is defined.

The {\em Dirichlet product} $f\otimes g$ will be defined by the
development
\begin{equation}\label{Dirprod} f\otimes g \;\; = \;\;
\sum_{\alpha} d_{\alpha}\phi_{\alpha} \;\; =
\;\;
\sum_{\alpha}\left(\sum_{\alpha_{1}\alpha_{2}=\alpha }a_{\alpha_{1}}b_{\alpha_{2}}\right) \phi_{\alpha},
\end{equation}
provided that it converges.  Let us say a few words about what this means.  Here, we note that for $\alpha \not=0$, we have that
$d_{\alpha}=\langle f, \tau_{\alpha}\bar{g}\rangle$ where 
$\tau_{\alpha}\bar{g} = \sum_{\beta\in K^{\ast}} b_{\alpha\beta^{-1}}\phi_{\beta}$, hence
is well-defined.  Consider also the formal expression
\[   d_{0} = a_{0}\sum_{\alpha\in K^{\ast}} b_{\alpha}  + b_{0}\sum_{\alpha\in K^{\ast}} a_{\alpha} 
+a_{0}b_{0}.
 \]
 Then we say that the Dirichlet product $f\otimes g$ converges when $d_{0}$ converges absolutely and the sequence $\{ d_{\alpha}\}$ defines an element of $l^{2}(K )$.
We denote by ${\sf Dom}_{\otimes}(f)$ the set of functions having defined Dirichlet product
with $f$.

In terms of the $ {\boldsymbol \zeta}$ parameter:
\[ (f\otimes g)( {\boldsymbol \zeta})
\;\;=\;\;\sum_{\alpha} a_{\alpha}\cdot  g\big( {\boldsymbol \zeta}^{\alpha}\big)\;\;=\;\;
\sum_{\alpha} b_{\alpha} \cdot f\big( {\boldsymbol \zeta}^{\alpha}\big) ,\]
showing that Dirichlet multiplication is commutative and distributive over
ordinary addition $+$ of functions.  It follows that $L^{2}(\SO_{K},\C )$
is a {\it partial algebra} with respect to each
of the operations $\oplus$ and $\otimes$ separately: here, {\it partial} refers
to the fact that  $\oplus$ and $\otimes$ are only partially defined,
and when they are, the usual axioms of group algebra hold for each operation.
When $\mathcal{K}$ is of infinite degree over $\Q$, 
the relevant discussion
applies to the space of square integrable functions
on the proto solenoid $L^{2}(\hat{\SO}_{\mathcal{K}},\C)$.

\begin{note}\label{mellintransform}  Suppose that $K=\Q$ and $f,g\in C_{0}(\T_{\Q},\C )\subset
C_{0}(\SO_{\Q},\C )$ with $a_{n}=b_{n}=0$ for $n\leq 0$.  Then
$f$ and $g$ define via a Mellin-type transform convergent Dirichlet series
\[  {\sf D}_{f}(y)\;\;=\;\;\sum a_{n}n^{-2\pi iy}\quad\quad\quad\quad
{\sf D}_{g}(y)\;\;=\;\;\sum b_{n}n^{-2\pi iy} ,  \]
$y\in\R$, and 
\[ D_{f\otimes g}\;\; =\;\; D_{f}\oplus D_{g}. \]
Thus the algebra of zeta functions, L-functions, Dirichlet series
{\em etc.}\ is codified by the Dirichlet product.  
\end{note}

Now let $\HSO_{K}$ be the hyperbolized adele class group defined in \S 6. 
$\HSO_{K}$ is a lamination whose leaves are $d$-dimensional polydisks,
and so we say that a continuous function $F:\HSO_{K}\rightarrow\C$
is holomorphic if its restriction $F|L$ to each leaf $L$ is
holomorphic, or equivalently (since all leaves are dense),
if its restriction to the canonical leaf $\HP_{K}=\HP_{K}^{\R}\times \HP_{K}^{\C}$ is holomorphic.  We recall that holomorphicity in
the factor $\HP_{K}^{\C}$ is defined in terms of the ``upper half-plane $\times$ right half-plane" multi-variable ${\bf u}\times{\bf v}$.

For each ${\bf t}\in  (0,\infty )^{d}$ let 
$\SO_{K}({\bf t})\subset\HSO_{K}$ be the subspace of points having 
extended coordinate
${\bf t}$ with respect to the decomposition (\ref{hyperbolized}).  Since $\SO_{K}({\bf t})\approx\SO_{K}$ we may put
on $\SO_{K}({\bf t})$ the unit mass Haar measure and define 
for $F,G :\HSO_{K}\rightarrow \C$ the pairing
\[ (F,G)_{\bf t}\;\;=\;\ \int_{\SO_{K}({\bf t})} 
F\,\overline{G}\; d\mu . \]

\begin{defi}  The {\bf Hardy space} associated to $K/\Q$ 
is the Hilbert space
\[ {\sf H}[K]\;\;=\;\; \Big\{ F:\HSO_{K}\rightarrow\C\;\Big| 
\;\;F\text{ is holomorphic and } \sup_{\bf t} (F,F)_{\bf t}<\infty\Big\}. \]
with inner-product 
$(\cdot ,\cdot )=\sup_{\bf t}(\cdot ,\cdot )_{\bf t}$.  
\end{defi}
Evidently any 
$F\in {\sf H}[K]$
has an {\em a.e.}\ defined $L^{2}$ limit for
${\bf t}\rightarrow {\bf 0}$, and such
a limit defines an element
\[ \partial F := f\in L^{2}(\SO_{K},\C ).\]
Using the Fourier development
available there, we may write the restriction $F|_{\HP_{K}}$ as follows: if ${\boldsymbol \rho}={\boldsymbol \tau}\times {\boldsymbol \kappa}\in \HP_{K}=\HP_{K}^{\R}\times \HP^{\C}_{K}$
(see \S \ref{hyperbolizations} for the relevant notation) then
\begin{equation}\label{FourRep} F|_{\HP_{K}}({\boldsymbol \rho})\;\;=\;\;\sum a_{\alpha}
\exp (2\pi i {\sf Tr}(\alpha\cdot {\boldsymbol \rho} ))\;\;=\;\;\sum a_{\alpha}
\exp (2\pi i {\sf Tr}(\alpha\cdot {\boldsymbol \tau} ))\exp (2\pi i {\sf Tr}(\alpha\cdot {\boldsymbol \kappa} ))
\end{equation}
where
\begin{equation}\label{FourRep1}
 \exp (2\pi i {\sf Tr}(\alpha\cdot {\boldsymbol \tau})) = \prod_{\nu} \exp(2\pi i \alpha_{\nu}x_{\nu} )\exp (-2\pi \alpha_{\nu}t_{\nu})
 \end{equation}
and 
\begin{equation}\label{FourRep2}
\exp (2\pi i {\sf Tr}(\alpha\cdot {\boldsymbol \kappa}  ))=
\prod_{(\mu , \bar{\mu})}  \exp(4\pi i ({\rm Re}(\alpha_{\mu}z_{\mu} ))
\exp (-4\pi {\rm Im}(\alpha_{\mu}b_{\mu} ))  .
\end{equation}   
We will sometimes switch to power series notation and shorten this to
\[ F|_{\HP_{K}}( {\boldsymbol \varrho} )\;\;=\;\; \sum_{\alpha}a_{\alpha} {\boldsymbol \varrho} ^{\alpha} \;=\;\; \sum_{\alpha}a_{\alpha} {\boldsymbol \xi}^{\alpha} {\boldsymbol \eta}^{\alpha},
\]
where
$ {\boldsymbol \varrho}\equiv \exp (2\pi i{\sf Tr}( {\boldsymbol \rho}) )$, $\ {\boldsymbol \xi}\equiv \exp (2\pi i{\sf Tr}( {\boldsymbol \tau} ))$ and 
$ {\boldsymbol \eta}\equiv \exp (2\pi i{\sf Tr}( {\boldsymbol \kappa}) )$.

By the {\it positive cone} in $K_{\infty}$ we shall mean the set of ${\bf x} = (x_{\nu};(z_{\mu}, \bar{z}_{\mu}))\in K_{\infty}$ for which $x_{\nu}> 0$
for all $\nu$ real and  $(z_{\mu}, \bar{z}_{\mu})\in\B\times\bar{\B}$ for each complex place
pair $(\mu ,\bar{\mu})$.  
It is clear then that for the series of (\ref{FourRep}) to define elements of ${\sf H}[K]$,
we must have
that $a_{\alpha}=0$ for $\alpha$ not contained in the positive
cone of $K_{\infty}$.
We have that 
$\| F\|^{2}=\sum |a_{\alpha}|^{2}$ and hence the correspondence $F\mapsto \partial F$ yields an
isometric inclusion of Hilbert spaces
\[ {\sf H}[K]\;\hookrightarrow\; L^{2}(\SO_{K},\C ) . \]

\begin{note}\label{singularcollapse} In (\ref{FourRep2}), observe that when $\alpha$ is real and positive (i.e. all of its place coordinates are real and positive),
then
$ {\boldsymbol \eta}^{\alpha}$ is a holomorphic function purely of the upper half plane
variable ${\bf u}=(u_{\mu})\in
\prod \HP_{\mu}$ of $\HP_{K}$ {\it i.e.}\
 is constant in the right plane variable ${\bf v}=(v_{\mu})\in\prod(-i\HP_{\mu})$.  Similarly, when $\alpha$ belongs to the positive
 imaginary axis,
 $ {\boldsymbol \eta}^{\alpha}$ only depends on ${\bf v}$.  In particular, this shows that
 $F|_{\HP_{K}^{\C}}$ is holomorphic with respect to the multi-variable ${\bf u}\times {\bf v}\in\HP_{K}^{\C}$.  
 \end{note}

In order to build up from ${\sf H}[K]$ a kind of holomorphic extension
of $K$, 
it will be necessary for us to be able to 
interpret all elements of $L^{2}(\SO_{K},\C )$ 
as boundaries of holomorphic functions.  We 
shall thus expand upon the $\Z/2\Z$-graded (super) 
convention, which has the virtue of regarding 
holomorphic and anti-holomorphic
functions on equal footing.  

\subsection{The Totally Real Case}\label{GradedSectReal}

We first consider the case $K=\Q$.  
Let $\Theta=\{ -,+\}\cong \Z/2\Z$ be the sign group and denote 
by $\HP^{-}=\R\times (-\infty ,0)$ the hyperbolic lower half-plane,
${\sf c}_{-}:\HP\rightarrow 
\HP^{-}$ complex-conjugation. Then every element $f=\sum a_{q}\zeta^{q}\in L^{2}(\SO_{\Q}, \C )$
determines a triple $(F_{-},F_{0}, F_{+})$
\[ F_{+}(\tau )\;\; =\;\; \sum_{q>0} a_{q}\exp (2\pi iq \tau),\quad
\quad  F_{0}\;\;=\;\;  a_{0},\quad\quad 
F_{-}(\tau)\;\; =\;\; \sum_{q<0} a_{q}\exp (2\pi iq {\sf c}_{-}(\tau)).
\]
The functions $F_{+}(z)$ and $F_{-}(z)$ are viewed as elements
of the Hardy spaces ${\sf H}_{+}[ \Q]= {\sf H}[\Q ]$ and 
${\sf H}_{-}[\Q ]$ = Hardy space of anti-holomorphic functions
on $\HSO_{\Q}$.

Now let $K/\Q$ be a totally real extension of degree $d$.  Let
\[\Theta_{K}=\{ -, +\}^{d}\cong (\Z/2\Z )^{d}\]
and write $\C_{K}=K_{\infty}\otimes_{\R} \C\cong \C^{d}$, whose points
are written ${\boldsymbol \tau}= (x_{\nu_{1}}+it_{\nu_{1}},\dots ,x_{\nu_{d}}+it_{\nu_{d}})$.
For each ${\boldsymbol \theta}\in
\Theta_{K}$,  define
\[ \HP_{K}^{{\boldsymbol \theta}}\;\;=\;\;\big\{ {\boldsymbol \tau}\in \C_{K}
\;\big|\;\;
\big(\text{sign}(t_{\nu_{1}}),\dots , \text{sign}(t_{\nu_{d}})) = {\boldsymbol \theta}
\big\} .\] 
We associate  to ${\boldsymbol \theta}$ a conjugation
${\sf c}_{{\boldsymbol \theta}}:\C_{K}\rightarrow\C_{K}$, where the coordinates
of  ${\sf c}_{{\boldsymbol \theta}}({\boldsymbol \tau})={\boldsymbol \tau}'$
satisfy
\[ x'_{\nu_{j}} +it'_{\nu_{j}} =x_{\nu_{j}} + \theta_{j}it_{\nu_{j}} \]
for $j=1,\dots ,d$.  The conjugation maps are not holomorphic in ${\boldsymbol \tau}$ and satisfy
for all ${\boldsymbol \theta}, {\boldsymbol \theta}_{1},{\boldsymbol \theta}_{2}\in\Theta_{K}$:
\[  {\sf c}_{{\boldsymbol \theta}}(\HP_{K}) = \HP_{K}^{{\boldsymbol \theta}} \quad
\text{and}\quad
{\sf c}_{{\boldsymbol \theta}_{1}}\circ{\sf c}_{{\boldsymbol \theta}_{2}}
={\sf c}_{{\boldsymbol \theta}_{1}{\boldsymbol \theta}_{2}}.\]
A   {\em $\theta$-holomorphic function} is one of the form $F\circ c_{\theta}$, where $F: \C_{K}\rightarrow \C$ is holomorphic.

Denote 
by $K^{{\boldsymbol \theta}}$ those elements $\alpha\in K$ whose coordinates with
respect to the embedding $K\hookrightarrow K_{\infty}$ satisfy 
\[ \big(\text{sign}(\alpha_{\nu_{1}}),\dots , \text{sign}(\alpha_{\nu_{d}})\big) \;\; =\;\; {\boldsymbol \theta} \]
Note that  
\[ K^{{\boldsymbol \theta}_{1}}\cdot K^{{\boldsymbol \theta}_{2}}\subset 
K^{{\boldsymbol \theta}_{1}{\boldsymbol \theta}_{2}}.\]   

Then every element $f=\sum a_{\alpha}{\boldsymbol \zeta}^{\alpha}\in L^{2}\in L^{2}(\SO_{K},\C )$ determines
a $(2^{d}+1)$-tuple $(F_{{\boldsymbol \theta}}; F_{0})$, where for each 
${\boldsymbol \theta}\in\Theta_{K}$, $F_{{\boldsymbol \theta}}:\HSO_{K}\rightarrow\C$ is defined
as the unique extension to $\HSO_{K}$ of the following ${\boldsymbol \theta}$-holomorphic function on $\HP_{K}$:
\[ F_{{\boldsymbol \theta}}({\boldsymbol \tau})\;\; =\;\; \sum_{\alpha\in K^{{\boldsymbol \theta}}} 
a_{\alpha} \exp \big(2\pi i  {\rm Tr}( \alpha \cdot {\sf c}_{{\boldsymbol \theta}}({\boldsymbol \tau})) \big)
\;\;\equiv\;\;\sum_{\alpha\in K^{{\boldsymbol \theta}}}a_{\alpha}({\sf c}_{{\boldsymbol \theta}}({\boldsymbol \xi}))^{\alpha} \]
where ${\sf c}_{{\boldsymbol \theta}}({\boldsymbol \xi})\equiv 
\exp (2\pi i {\sf Tr}({\sf c}_{{\boldsymbol \theta}}({\boldsymbol \tau})))$.The term 
$F_{0}$ is the constant function $a_{0}$. 

The Hardy space of $\theta$-holomorphic functions is 
denoted ${\sf H}_{{\boldsymbol \theta}}[K]$.
The space of $(2^{d}+1)$-tuples
is viewed as a graded Hilbert space:
\[ {\sf H}_{\bullet}[K]\;\; = 
\;\;\Big(\bigoplus_{{\boldsymbol \theta}}{\sf H}_{{\boldsymbol \theta}}[K]\Big)\oplus\C ,
\]
whose inner-product is the direct sum of the inner-products
on each of the summands.
We will often write ${\sf H}[K]$ for the summand of
${\sf H}_{\bullet}[K]$ corresponding to ${\boldsymbol \theta} =(+,\dots ,+)$
{\em i.e.}\ the Hardy space of functions holomorphic on $\HSO_{K}$
in the ordinary sense.

The Cauchy and Dirichlet products are defined on ${\sf H}_{\bullet}[K]$
via boundary extensions {\em e.g.}\ $F\otimes G$ is defined
to be the unique element of ${\sf H}_{\bullet}[K]$ whose
boundary is $\partial F\otimes\partial G$, provided the latter is defined.  
The Cauchy product does not generally
respect the grading, but the Dirichlet product 
has the following graded decomposition law:
\begin{equation}\label{gradepreserve} \big(F\otimes G\big)_{{\boldsymbol \theta}} \;\;=\;\;
\sum_{{\boldsymbol \theta}={\boldsymbol \theta}_{1}{\boldsymbol \theta}_{2}}F_{{\boldsymbol \theta}_{1}}\otimes G_{{\boldsymbol \theta}_{2}} 
\quad\quad\quad\quad
\big( F\otimes G\big)_{0} \;\;=\;\;
 F(1) G(1)\;\; -\;\; F_{0} G_{0}.
 \end{equation}
 
 \subsection{The Totally Complex Case}\label{GradedSectComplex}

Now let us suppose that $K/\Q$ is totally complex of degree $d=2s$.  In order to extend the ideas in
the previous paragraphs in a way compatible with the complex places,
we develop a complex theory of signs.  We consider the {\it singular sign group}
\[\Theta^{\C}= \{ \sqrt{-}, -, -\sqrt{-}, +\} \cong \Z/4\Z\] 
and say that $z\in(\R\cup i\R)-0$ has singular sign $\sqrt{-}, -, -\sqrt{-}, +$ according to whether
$z$ belongs to $i\R_{+}$, $ -\R_{+}$, $-i\R_{+}$ or $\R_{+}$.  

For points which do not belong to $\R\cup i\R$, we introduce the {\it complex sign set}
\[ \Omega = \{  \sqrt{-}\varepsilon, -\varepsilon, -\sqrt{-}\varepsilon, +\varepsilon   \},\]
viewed as a $\Theta^{\C}$-set in the obvious way.
We say that $z\in\C - (\R\cup i\R)$ has complex sign 
$\sqrt{-}\varepsilon, -\varepsilon, -\sqrt{-}\varepsilon, +\varepsilon$ according to whether
$z$ belongs to $i\B$, $-\B$, $-i\B$ or $\B$, where $\B$ is as before the quarter plane
of $z=x+iy$ with $x,y>0$.  
Every $\omega\in \Omega$ can be written uniquely in the form
\[  \omega = \theta \varepsilon  \]
for $\theta\in  \Theta^{\C}$, and the map
\[   e:\Omega \rightarrow  \Theta^{\C},\quad e(\theta\varepsilon ) =\theta \]
 is an isomorphism of $ \Theta^{\C}$-sets.   The singular sign group $\Theta^{\C}$ may be viewed intermediate to
the ``signless" element $0$ and the complex sign $\Theta^{\C}$-set  $\Omega$.

 If $z,z'$ have complex signs $\omega, \omega '$ then the product can have sign $e(\omega )e(\omega')\varepsilon$,
or else it can ``overflow" into the neighboring signs: the singular sign
$\sqrt{-}e(\omega )e(\omega')$ and the complex sign $\sqrt{-}e(\omega )e(\omega')\varepsilon$: explaining why we view $\Omega$ as no more
than a $\Theta^{\C}$-set .  We view the union 
\[ \mho := \Theta^{\C}\cup \Omega\] as an abelian ``multi-valued
group" i.e.\ a set equipped with the abelian set-valued product
\[ \mho\times\mho \longrightarrow {\sf 2}^{\mho}   \]
with identity the sign $+$.
We write for $\vartheta, \vartheta_{1},  \vartheta_{2}\in\mho$
\[  \vartheta\in \vartheta_{1}\cdot \vartheta_{2} \]
to indicate that $\vartheta$ is among the possible signs that the product $z_{1}z_{2}$
can assume when ${\sf sign}(z_{1})=\vartheta_{1}$, 
${\sf sign}(z_{2})=\vartheta_{2}$.  Note that $|\vartheta_{1}\cdot \vartheta_{2} |>1$
only when $\vartheta_{1}, \vartheta_{2}\in\Omega$, in which case 
$|\vartheta_{1}\cdot \vartheta_{2} |=3$.
The function $e$ of the previous paragraph extends to $ \mho$
by the identity in $\Theta^{\C}$. 
 
 Consider now the set 
 \[ \mho_{K}=\bigg\{ {\boldsymbol \vartheta} = (\vartheta_{j}) \in \mho^{s}|\; \exists \alpha\in K\text{ such that }   
    {\rm sign}(\alpha_{\mu_{j}} )= \vartheta_{j} \text{ for } j=1,\dots ,s \bigg\}.\]
  Notice that we are only using the first element $\mu_{j}$ of each complex place pair $(\mu_{j} ,\bar{\mu}_{j})$ 
  to define the sign, and that we have that for all $K$ the sign $(+,\cdots ,+)\in \mho_{K}$.  
  The set-valued product of components induces a set-valued product
  \[\mho_{K}\times \mho_{K}\longrightarrow {\sf 2}^{\mho_{K}}.  \]
    The map 
  \[  e:\mho_{K}\longrightarrow (\Theta^{\C})^{s}\]
  is defined $e({\boldsymbol \vartheta})=(e(\vartheta_{j}))$.

   If we define the {\it type} of ${\boldsymbol \vartheta}\in\mho_{K} $
  to be the vector $t({\boldsymbol \vartheta})$ whose $j$th component is $1$ or $\varepsilon$ depending on whether
  $\vartheta_{j}\in \Theta^{\C}$ or $ \Omega$, then
  \[ {\boldsymbol \vartheta} = e({\boldsymbol \vartheta})\cdot t({\boldsymbol \vartheta}) .\]
  We will denote ${\boldsymbol \varepsilon}= (\varepsilon ,\dots ,\varepsilon )$ and ${\bf 1}= (1,\dots ,1 )$, and say that
 ${\boldsymbol \vartheta}$ is {\it 
  complex homogeneous} ({\it 
  singular homogeneous}) if $t({\boldsymbol \vartheta}) = {\boldsymbol \varepsilon}$ ($t({\boldsymbol \vartheta})={\bf 1}$).  If ${\boldsymbol \vartheta}$ is
  singular homogeneous, we will sometimes write ${\boldsymbol \vartheta} = {\boldsymbol \theta}$.
  The singular homogeneous elements form a subgroup $\Theta^{\C}_{K}$ with respect to which 
  $\mho_{K}$ is a $\Theta^{\C}_{K}$-set.

  \begin{exam}  Let $K$ be the splitting field for $X^{3}-2$.  Then the sign of any root is (modulo ordering of the places)
  $(+,   \sqrt{-}\varepsilon , -  \varepsilon)$, so that $\mho_{K}$ has inhomogeneous triples.  On the other hand,
  if $K$ is the splitting field of $X^{3}-1$ then all elements of $K$ are homogeneous.
  This is also true of the Gaussian numbers $\Q (i)$, where $ \mho_{\Q (i)} = \mho$.
  \end{exam}

For each ${\boldsymbol \vartheta}\in \mho_{K}$ we will associate a product of half spaces.
For $\omega\in\Omega$ and each place pair $(\mu ,\bar{\mu})$ write
\[ \HP_{(\mu ,\bar{\mu})}^{\omega}\;\;=\;\;\big\{  (\kappa_{\mu},\bar{\kappa}_{\mu})\in \C^{2}\times\C^{2}
\;\big|\;\;
\text{sign}(b_{\mu}) = \omega
\big\} \]
where we recall that $\kappa_{\mu} =(z_{\mu}, b_{\mu})$, see \S \ref{hyperbolizations}.
Notice that when $\omega=+\epsilon$, $\HP_{(\mu ,\bar{\mu})}^{+\epsilon}=\HP_{(\mu ,\bar{\mu})}$.  We will view
these spaces as products of half-spaces in the same way that we viewed $\HP_{(\mu ,\bar{\mu})}$
as a product of half-spaces, see (\ref{upperright}) in \S \ref{hyperbolizations}.  For example, when $\omega = \sqrt{-}\varepsilon$,
$b_{\mu}=s_{\mu}+it_{\mu}$ satisfies $s_{\mu}<0$ and $t_{\mu}>0$ so that
we obtain the product of the upper half-plane with the left half-plane:
\[  \HP_{(\mu ,\bar{\mu})}^{\sqrt{-}\varepsilon} = \HP_{\mu}\oplus i\HP_{\mu}. \]
The points of the first factor are denoted $u_{\mu}=x_{\mu}+it_{\mu}$ and those of the second factor by $v_{\mu}=s_{\mu}+iy_{\mu}$.
In the same fashion, we have
\[   \HP_{(\mu ,\bar{\mu})}^{-\varepsilon} = -\HP_{\mu}\oplus i\HP_{\mu},\quad 
 \HP_{(\mu ,\bar{\mu})}^{-\sqrt{-}\varepsilon} = -\HP_{\mu}\oplus -i\HP_{\mu}. \]
Also, for each singular sign $\theta\in \Theta^{\C}$, we write 
\[ \HP_{(\mu ,\bar{\mu})}^{\theta} :=  \HP_{(\mu ,\bar{\mu})}^{\theta\varepsilon}. \]
Then for ${\boldsymbol \vartheta}\in \mho_{K}$ we define the ${\boldsymbol \vartheta}${\it -hyperbolic plane} as 
\[  \HP_{K}^{{\boldsymbol \vartheta}} =  \prod_{j=1}^{s}  \HP_{(\mu_{j} ,\bar{\mu}_{j})}^{\vartheta_{j}}. \]

We now define conjugation maps relating these planes.  For $\kappa_{\mu}=(z_{\mu},b_{\mu})\in\C^{2}$ and
 $\theta\in\Theta^{\C}$ we define
 \[  {\sf c}_{\theta} (z_{\mu},b_{\mu}) = (z_{\mu},\theta b_{\mu}) . \]
This induces a map of $\HP_{(\mu ,\bar{\mu})}=\{ (z,\bar{z})\times (b, -\bar{b} )\;|\; z\in\C, b\in\B \}$ such that
 \[  {\sf c}_{\theta}\big( \HP_{(\mu ,\bar{\mu})}\big) = \HP_{(\mu ,\bar{\mu})}^{\theta}= \HP_{(\mu ,\bar{\mu})}^{\theta\varepsilon} .\]
 Notice that the conjugation maps are not (holomorphic) functions of the half-plane 
 variables $u_{\mu}=x_{\mu}+it_{\mu}, v_{\mu}=s_{\mu}+iy_{\mu}$.
 For example, when $\theta = \sqrt{-}$, 
 \[   {\sf c}_{\sqrt{-}}  (x_{\mu}+it_{\mu}, s_{\mu}+iy_{\mu})= (x_{\mu} +is_{\mu},-t_{\mu}+iy_{\mu} ). \]
We have
 \[   {\sf c}_{\theta_{1}} \circ  {\sf c}_{\theta_{2}} = {\sf c}_{\theta_{1}\theta_{2}}\]
 for all $\theta_{1}, \theta_{2}\in\Theta^{\C}$.

For each ${\boldsymbol \vartheta}\in \mho_{K}$ we associate a conjugation
${\sf c}_{{\boldsymbol \vartheta} }:\C_{K}\rightarrow\C_{K}$
where the coordinates of 
${\sf c}_{{\boldsymbol \vartheta} }({\boldsymbol \kappa})= {\boldsymbol \kappa}' $ are determined by
\begin{equation}  \kappa'_{\mu_{j}}= (z'_{\mu_{j}}, b'_{\mu_{j}}) = (z_{\mu_{j}}, e(\vartheta_{j} ) b_{\mu_{j}}) \end{equation}
for $k=1,\dots ,s$.  Notice that ${\sf c}_{{\boldsymbol \vartheta} }$ restricts to a map
${\sf c}_{{\boldsymbol \vartheta} }:\HP_{K}\rightarrow\HP_{K}^{{\boldsymbol \vartheta}}$. 
A function $F:\C_{K}\rightarrow \C$ of the shape $F=G\circ {\sf c}_{{\boldsymbol \vartheta} }$
where $G:\C_{K}\rightarrow \C$ is holomorphic, will be called {\em ${\boldsymbol \vartheta}$-holomorphic}.

For each ${\boldsymbol \vartheta}\in\mho_{K}$, denote by $K^{{\boldsymbol \vartheta}}$ those elements $\alpha$ whose coordinates with
respect to the embedding $K\hookrightarrow K_{\infty}$ satisfy 
\[ \big(\text{sign}(\alpha_{\mu_{1}}),\dots , \text{sign}(\alpha_{\mu_{s}})\big) \;\; =\;\; {\boldsymbol \vartheta} .\]
For  ${\boldsymbol \theta}\in \Theta^{\C}_{K}$ and arbitrary ${\boldsymbol \vartheta}\in \mho_{K}$ we have
\begin{equation}\label{prodlaw}   K^{{\boldsymbol \theta}}\cdot K^{{\boldsymbol \vartheta}}\subset K^{{\boldsymbol \theta} {\boldsymbol \vartheta}} .
\end{equation}

In general, we cannot expect a law of the genre (\ref{prodlaw}) since the complex signs do not form a group.  We have instead the multi-valued law:
\[ K^{{\boldsymbol \vartheta}_{1}}\cdot K^{{\boldsymbol \vartheta}_{2}}\subset 
\bigcup_{{\boldsymbol \vartheta}\in {\boldsymbol \vartheta}_{1}\cdot {\boldsymbol \vartheta}_{2}  } K^{ {\boldsymbol \vartheta}} .  \]

Every element $f\in L^{2}(\SO_{K},\C )$ now determines
a $(|\mho_{K}|  +1)$-tuple $(F_{\boldsymbol \vartheta};F_{0})$ as follows.
For each 
${\boldsymbol \vartheta}\in\mho_{K}$, $F_{\boldsymbol \vartheta}:\HSO_{K}\rightarrow\C$ is defined
as the extension to $\HSO_{K}$ of the following function on $\HP_{K}$:
\[F_{\boldsymbol \vartheta}({\boldsymbol \kappa})\;\; =\;\; \sum_{\alpha\in K^{\boldsymbol \vartheta}} 
a_{\alpha} \exp \big(2\pi i  {\sf Tr}( \alpha {\sf c}_{\boldsymbol \vartheta}({\boldsymbol \kappa})) \big)
\;\;\equiv\;\;\sum_{\alpha\in K^{\boldsymbol \vartheta}}a_{\alpha}({\sf c}_{\boldsymbol \vartheta}({\boldsymbol \eta}))^{\alpha} .\]
where  ${\sf c}_{\boldsymbol \vartheta}({\boldsymbol \eta})\equiv 
\exp (2\pi i {\sf Tr}({\sf c}_{\boldsymbol \vartheta}({\boldsymbol \kappa})))$.
Observe that the comments found in {\it Note} \ref{singularcollapse} imply that whenever the sign coordinate $\vartheta_{j}$ is singular i.e.\  $\vartheta_{j}\in\Theta^{\C}$
then $F_{\boldsymbol \vartheta}$ is constant in one of the corresponding half-plane coordinates $u_{\mu_{j}}, v_{\mu_{j}}$.
This explains why we call the signs in $\Theta^{\C}$ ``singular''.

We refer to $F_{\boldsymbol \vartheta}$ as the
$\boldsymbol \vartheta$-holomorphic component of $f$,
and the Hardy space of $\boldsymbol \vartheta$-holomorphic  functions is 
denoted ${\sf H}_{\boldsymbol \vartheta}[K]$.
We obtain a graded Hilbert space:
\[ {\sf H}_{\bullet}[K]\;\; = 
\;\;\Big(\bigoplus_{\boldsymbol \vartheta}{\sf H}_{\boldsymbol \vartheta}[K]\Big)\oplus \C ,
\]
whose inner-product is the direct sum of the inner-products
on each of the summands.
We will write as in the real case ${\sf H}[K]$ for the summand of
${\sf H}_{\bullet}[K]$ corresponding to $\boldsymbol \vartheta ={\bf 1}$.

The Cauchy and Dirichlet products are defined on ${\sf H}_{\bullet}[K]$
via boundary extensions just as in the real case.
When ${\boldsymbol \vartheta}\in{\boldsymbol \vartheta}_{1}\cdot {\boldsymbol \vartheta}_{2}$ we write 
\[   (F_{{\boldsymbol \vartheta}_{1}}\otimes G_{{\boldsymbol \vartheta}_{2}} )|_{{\boldsymbol \vartheta}} \]
for the projection of $F_{{\boldsymbol \vartheta}_{1}}\otimes G_{{\boldsymbol \vartheta}_{2}} $ onto the sub series
indexed by $\alpha\in K^{{\boldsymbol \vartheta}}$.
Then we have the following graded decomposition law generalizing that 
described in (\ref{gradepreserve}) for the totally real case:
\begin{equation}\label{gradepreserve2} \big(F\otimes G\big)_{{\boldsymbol \vartheta}}  = 
\sum_{{\boldsymbol \vartheta}\in{\boldsymbol \vartheta}_{1}{\boldsymbol \vartheta}_{2}} \mathfrak(F_{{\boldsymbol \vartheta}_{1}}\otimes G_{{\boldsymbol \vartheta}_{2}})|_{{\boldsymbol \vartheta}},\quad  \big( F\otimes G\big)_{0} \;\;=\;\;
 F(1) G(1)\;\; -\;\; F_{0} G_{0}. \end{equation}
 
 \vspace{8mm}
 
Now consider $K/\Q$ a general finite extension.  The we obtain a decomposition of the form
 \begin{equation}\label{realcomplex}  {\sf H}_{\bullet}[K] = {\sf H}^{\R}_{\bullet}[K]\oplus {\sf H}^{\C}_{\bullet}[K]\end{equation}
corresponding to the real and complex places (graded accordingly) so that in particular, every 
$f\in L^{2}(\SO_{K},\C )$  determines
a $(2^{r}+ |\mho_{K}|  +1)$-tuple of functions.

\begin{note}  Let $K/\Q$ be Galois.  Then the action of the Galois group
${\sf Gal}(K/\Q )$ induces a well-defined action on the sign group $\Theta_{K}$ (when
$K$ is real) and on the ``multi-valued sign group" $\mho_{K}$ 
(when $K$ is complex) .
\end{note}

For $K$ real or complex, we have a sub partial double group algebra
${\sf H}_{\bullet}[O_{K}]$ of $ {\sf H}_{\bullet}[K]$ defined by Fourier series whose 
indices belong only
to $O_{K}$.   The Hilbert space of graded holomorphic functions 
pulled back from the Minkowski hyperbolized torus $\mathfrak{T}_{K}$
are denoted
\[    {\sf H}_{\bullet}[\mathfrak{T}_{K}] :={\sf H}_{\bullet}[\mathfrak{d}_{K}^{-1}]  , \]
where $\mathfrak{d}_{K}^{-1}$ is the inverse different.  This sub Hilbert space is closed
with respect to the Cauchy product (where it is defined), and is closed with respect
to the action
by Dirichlet multiplication with elements of ${\sf H}_{\bullet}[O_{K}]$ (whenever such products
are defined).
When $K=\Q$,  then we have  
${\sf H}_{\bullet}[\Z]=  {\sf H}_{\bullet}[\mathfrak{T}_{\Q}] $ is
a partial double algebra with respect to the Cauchy and Dirichlet products.

We now indicate how to extend this construction to infinite field extensions $\mathcal{K}/\Q$.
Here we use the hyperbolized adele class group $\HSO_{\mathcal{K}}$ 
(associated to $\SO_{\mathcal{K}}$) in conjunction
with the proto adele class group $\hat{\SO}_{\mathcal{K}}$.  
A continuous function $F:\HSO_{\mathcal{K}}\rightarrow\C$
is holomorphic if its restriction to any of the dense leaves is holomorphic.
In particular, we note that the restriction $F|_{\HP_{\mathcal{K}}}$
is holomorphic if $F$ is holomorphic separately in each factor of the
polydisk decomposition $\HP_{\mathcal{K}}\approx \prod \HP_{\nu}\times\prod  \HP_{(\mu ,\bar{\mu})}$.
As in the case of a finite field extension, we define 
$\SO_{\mathcal{K}} ({\bf t})$ as the subset of $\HSO_{\mathcal{K}} $
having extended coordinate ${\bf t}\in (0,\infty )^{\infty}$.  The
proto compactification of $\SO_{\mathcal{K}} ({\bf t})$ is a lamination
$\hat{\SO}_{\mathcal{K}} ({\bf t})$
homeomorphic to $\hat{\SO}_{\mathcal{K}}$.  If we let
\[ \hat{\HSO}_{\mathcal{K}}\;\;=\;\;\hat{\SO}_{\mathcal{K}}\times
(0,\infty )^{\infty}\;\;=\;\;\bigcup\; \hat{\SO}_{\mathcal{K}} ({\bf t}),\]
the Hardy
space ${\sf H}[\mathcal{K}]$ is defined as the space
of holomorphic functions $F:\HSO_{\mathcal{K}}\rightarrow\C$
having a continuous extension 
$\hat{F}:\hat{\HSO}_{\mathcal{K}} \rightarrow\C$
and for which the norm
\[ \| \hat{F}\|_{\bf t}^{2}\;\;=\;\;
\int_{\hat{\SO}_{\mathcal{K}}({\bf t}) } |\hat{F}|^{2}\,d\mu \] 
is uniformly bounded in ${\bf t}$ (where $d\mu$
is induced from the Haar measure $\mu$ on 
$\hat{\SO}_{\mathcal{K}}\approx\hat{\SO}_{\mathcal{K}}({\bf t})$).
As in the finite-dimensional case, ${\sf H}[\mathcal{K}]$
is a Hilbert space with respect to the supremum of the integration
pairings on each $\hat{\SO}_{\mathcal{K}} ({\bf t})$.
The rest of the development follows that of the finite
extension case, where the grading is defined for the real and complex
places separately.

We summarize the above remarks in the following

\begin{theo}  Let $\mathcal{K}$ be a (possibly infinite degree) 
algebraic number field over $\Q$.  Then ${\sf H}_{\bullet}[\mathcal{K}]$
is a graded Hilbert space equipped with the structure of
a partial double $\C$-algebra with respect to the operations of $\oplus$ and
$\otimes$. 
\end{theo} 

 \section{Nonlinear Number Fields}\label{NonlinearChap}

Let $K/\Q$ be a number field of finite degree over $\Q$.  Let $\C [K]$ denote the vector space of formal, finite
$\C$-linear combinations of elements in $K$  i.e.\  
expressions
of the form $\sum a_{\alpha}\alpha$ where $\alpha\in K$ and $a_{\alpha}\in\C$, zero
for all but finitely many $\alpha$.
The operations $+_{K}$ and $\times_{K}$
extend linearly to $\C [K]$ yielding 
two operations, written $\oplus$ and $\otimes$, which
define on $\C [K]$
two algebra structures.  We refer to $\C [K]$ as the {\it field algebra}
generated by $K$.
To avoid confusion, we use the notation ${\sf id}_{\oplus}=0_{K}$
and ${\sf id}_{\otimes}=1_{K}$; $0$ will denote the vector space
identity, the element of $\C [K]$
for which $a_{\alpha}=0$ for all $\alpha$.
There exists a canonical double algebra monomorphism generated by
$\alpha\mapsto   {\boldsymbol \zeta}^{\alpha}$,
\[  \C [K]\;\hookrightarrow\;L^{2}(\SO_{K},\C )\cong {\sf H}_{\bullet} [K]  \]
having dense image.  
The subspace $\C [O_{K}]$ 
of $\C$-linear combinations of integers is
closed with respect to both algebra structures.   
These algebra structures are not compatible in any familiar
sense as the operations $\oplus$ and $\otimes$ {\em do not} 
obey the distributive law.

We define a linear map ${\sf T}:\C [K]\rightarrow\C $ by   
\[ {\sf T}(f)\;\;=\;\;\sum a_{\alpha}\;\in\;
\C.\]
Notice that
$f\otimes {\sf id}_{\oplus} = {\sf T}(f)	\cdot {\sf id}_{\oplus}$. 
The vector space $I_{K}={\sf Ker}({\sf T})$ is an ideal in $\C [K]$ 
with respect to the operations
of $\oplus$ and $\otimes$ and the set-theoretic difference
\[  \C^{\ast} [K]\;\; :=\;\;\C [K] - I_{K}\]
is preserved by both of $\oplus$ and $\otimes$.
 Denote a typical element of $\C^{\ast} [K]$ by $f^{\ast}$. 
Then for any $f^{\ast}$ we have
 \begin{equation}\label{nozerodiv}  f^{\ast}\otimes {\sf id}_{\oplus}^{\ast} 
 \;\;\in\;\; 
 \C^{\ast}\cdot {\sf id}^{\ast}_{\oplus}. \end{equation} 

Let ${\sf N}^{0} [K]$ denote the image of $\C^{\ast}[K ]$ in the complex
projectivization $\PR\C [K]$ of $\C [K]$. By virtue
of (\ref{nozerodiv}), the
operations $\oplus$ and $\otimes$ descend to ${\sf N}^{0} [K]$,
making it a {\em double semigroup}: a set with two
semigroup structures having no {\em a priori} compatibility.  We denote
its elements $[f]$.  The sub double semigroup ${\sf N}^{0} [O_{K}]$ is defined
similarly.

Note that the element $[{\sf id}_{\oplus}]$ behaves very much like the zero in a field
in that it is a universal annihilator with respect to $\otimes$: for all $[f]\in {\sf N}^{0} [K]$,
\[  [f]\otimes [{\sf id}_{\oplus}]\;\; =\;\; [{\sf id}_{\oplus}] .\]
Furthermore, the natural inclusions $O_{K}\hookrightarrow \C [O_{K}]$ and
$K\hookrightarrow \C [K]$ induce monomorphisms $O_{K}\hookrightarrow {\sf N}^{0} [O_{K}]$
and $K\hookrightarrow {\sf N}^{0} [K]$.  These echos with field theory make 
the double semi-group ${\sf N}^{0} [K]$ a natural paradigm for the 
ensuing development of nonlinear fields.
We remark that the preceding comments hold without change for
an infinite algebraic extension $\mathcal{K}/\Q$.
 
 Motivated by the monomorphism 
 $K\subset \C [K]\hookrightarrow {\sf H}_{\bullet}[K]$ 
 we set out to create from ${\sf H}_{\bullet}[K]$
 something akin to a 
 field extension of $K$ by graded holomorphic functions.  In this connection, we note that 
 the operations $\oplus$ and $\otimes$
 restrict to $+_{K}$ and $\times_{K}$ on $K$, and on the other hand,
 the vector space operations of point-wise addition and scalar multiplication
 do not preserve $K$.  Accordingly, we discard the
 vector space structure by projectivizing, retracing
 the steps in the construction of ${\sf N}^{0} [K]$.  
 
 Here, the trace operator ${\sf T}$ is not defined on all of 
 ${\sf H}_{\bullet}[K]$.  It is unambiguously defined on the subspace
 of functions $F$ having boundary $\partial F$ lying in the subspace $l^{1}(K)\cap l^{2}(K)$ of functions whose Fourier coefficients $\{ a_{\alpha}\}$
are absolutely summable.
 For such elements $F\in {\sf H}_{\bullet}[K]$ we define 
 ${\sf T}(F)=\sum a_{\alpha}=F(1)$ and
 denote 
 ${\sf I}_{K}={\sf Ker}({\sf T})$.
 We note that ${\sf I}_{K}$ is not closed in 
in ${\sf H}_{\bullet}[K]$ and is in fact dense there.
 The 
set theoretic difference
 \[ {\sf H}_{\bullet}^{\ast}[K]\;\; :=\;\;{\sf H}_{\bullet}[K]-{\sf I}_{K}
 \]
 inherits the grading
by restriction.  The associated quotient by $\C^{\ast}$, denoted 
\[ {\sf N}_{\bullet}[K],\] 
is an infinite dimensional subspace of the full projectivization $\PR{\sf H}_{\bullet}[K]$
which inherits $\oplus$, 
$\otimes$ as partially defined operations.  The grading gives rise to an arrangement of subspaces
\[ \{  {\sf N}_{\boldsymbol \vartheta}[K] \}  \]
where ${\sf N}_{\boldsymbol \vartheta}[K] = \PR{\sf H}_{\boldsymbol \vartheta}[K]\cap {\sf N}_{\bullet}[K]$.
The canonical monomorphism
\[  {\sf N}^{0} [K]\;\hookrightarrow\; {\sf N}_{\bullet}[K]\]
induced by $\alpha \mapsto {\boldsymbol \varrho}^{\alpha}$, has dense image,
and the elements of ${\sf N}^{0} [K]$ may be Cauchy or Dirichlet multiplied
with any element of ${\sf N}_{\bullet}[K]$.
We also have a monomorphism $K\hookrightarrow {\sf N}_{\bullet}[K]$.
Elements of  ${\sf N}_{\bullet}[K]$ will be denoted
by $[F]$.  
The sub partial double semigroup ${\sf N}_{\bullet}[O_{K}]$ 
is defined similarly. 
 These remarks are valid without change for an
infinite field extension $\mathcal{K}/\Q$.

 The grading of ${\sf N}_{\bullet}[K]$ induces one on ${\sf N}^{0} [K]$ and so we write
${\sf N}^{0}_{\bullet} [K]$.
 
 \begin{defi}\label{nonlinear}  A {\bf nonlinear number field} is a graded topological abelian partial double
 semigroup ${\sf S}_{\bullet}$ with respect to two operations $\oplus$ and $\otimes$ such that
 \begin{enumerate}
 \item  There exists a (possibly infinite degree) algebraic
 number field $\mathcal{K}/\Q$ and a graded double semigroup monomorphism 
 $\imath :{\sf N}^{0}_{\bullet} [\mathcal{K}]\hookrightarrow {\sf S}_{\bullet}$ having dense image.
 \item  The identity ${\sf id}_{\oplus}$ is a universal annhilator
 for $\otimes$: for all $F\in {\sf Dom}_{\otimes}({\sf id}_{\oplus} )$,
 $F\otimes {\sf id}_{\oplus}=  {\sf id}_{\oplus}$.  
 \end{enumerate}
 The closure ${\sf O}$ of the image $\imath ({\sf N}^{0} [O_{\mathcal{K}}])$ is called
 the {\bf nonlinear ring of integers}.
 \end{defi}
 
The qualificative ``nonlinear'' refers to the fact that distributivity in an ordinary
 field 
 is equivalent to the fact that the multiplication map is a bilinear
 operation. For $\mathcal{K}/\Q$ be a possibly infinite degree extension,
 notice that both ${\sf N}^{0}_{\bullet}[\mathcal{K}]$ and ${\sf N}_{\bullet}[\mathcal{K}]$
 are nonlinear number fields.  In addition the following are also nonlinear number fields:
  \begin{itemize}
 \item[-] $\bar{\sf N}_{\bullet}[\mathcal{K}]= \PR {\sf H}_{\bullet}[\mathcal{K}]$ = the full projectivization of
 ${\sf H}_{\bullet}[\mathcal{K}]$.  
 \item[-] Let ${\sf W}_{\bullet}[\mathcal{K}]$ be the projectivization of the subspace
 of ${\sf H}_{\bullet}[\mathcal{K}]$ consisting of absolutely convergent series with non-zero trace.  Then 
 ${\sf W}_{\bullet}[\mathcal{K}]$ is the {\it Wienerian nonlinear number field} associated to $\mathcal{K}$: a full (and not partial) subsemigroup of ${\sf N}_{\bullet}[\mathcal{K}]$.
 \end{itemize}
 We have inclusions $\mathcal{K}\subset {\sf N}^{0}_{\bullet}[\mathcal{K}]\subset{\sf N}_{\bullet}[\mathcal{K}]\subset
 {\sf W}_{\bullet}[\mathcal{K}]\subset \bar{\sf N}_{\bullet}[\mathcal{K}]$, the last three of which
 are dense.
 
 \begin{note}  In view of Note~\ref{mellintransform}, all of the arithmetic
 of classical (single variable) 
 zeta functions, Dirichlet series, $L$-functions, {\em etc.}\ is contained in
 ${\sf N}_{\bullet}[\mathcal{\Z}]$. 
 \end{note}

 \begin{note}  The set of Cauchy units
 ${\sf U}_{\bullet}[K]_{\oplus}$ form a dense subset of
 ${\sf N}_{\bullet}[K]$, since it contains all classes represented by 
real analytic nonvanishing
 functions of $\SO_{K}$.  It is an interesting question as to whether
 the set of Dirichlet units ${\sf U}_{\bullet}[K]_{\otimes}$
 is also dense, and whether there exists an algorithmic procedure
 to determine the coefficients of a Dirichlet inverse analogous to
 the classical M\"{o}bius inversion formula.
 \end{note}

 The following theorem says that on a dense subset,
 ${\sf N}_{\bullet}[\mathcal{K}]$ is the ``nonlinear field of fractions''
 of ${\sf N}_{\bullet}[O_{\mathcal{K}}]$:

\begin{theo}\label{fractions}  Let $\mathcal{K}$
 be a (possibly infinite degree) number field.  Then there is a dense subset
${\sf P}\subset{\sf N}_{\bullet}[\mathcal{K}]$ such that for all $[F]\in {\sf P}$,
there exists $[A]\in {\sf N}_{\bullet}[O_{\mathcal{K}}]$
with
\[ [A]\otimes [F]\;\in \; {\sf N}_{\bullet}[O_{\mathcal{K}}] .\] 
\end{theo}

\begin{proof}  Let ${\sf N}_{\bullet}[\mathcal{K}]_{\rm fin}$ be the subspace
associated to functions whose nonzero Fourier coefficients are
indexed by $\alpha$ in some fractional ideal $O_{\mathcal{K}}
\mathfrak{a}^{-1}\subset \mathcal{K}$, where $\mathfrak{a}\subset O_{\mathcal{K}}$.
Consider the sub double semigroup ${\sf N}_{\bullet}[\mathfrak{a}]$.
Then given $[F]\in {\sf N}_{\bullet}[\mathcal{K}]_{\rm fin}$ whose
nonzero Fourier coefficients are indexed by $O_{\mathcal{K}}\mathfrak{a}^{-1}$, there
exists $[A]\in {\sf N}_{\bullet}[\mathfrak{a}]$
such that $[A]\otimes [F]\in  {\sf N}_{\bullet}[O_{\mathcal{K}}]$.
\end{proof}
 
\section{Galois Groups and the Action of the Idele Class Group of $\Q$}\label{GalGroupChap}

Let $\mathcal{K}/\Q$ be a possibly infinite degree algebraic number field Galois over $\Q$.  In
this case, $\mathcal{K}$ is either totally real or totally complex: in either case we denote
by $\mho_{\mathcal{K}}$ the sign set.
We will also not distinguish $\mathcal{K}$ from its images in ${\sf H}_{\bullet}[\mathcal{K}]$ or
${\sf N}_{\bullet}[\mathcal{K}]$.  
We will prefer here to represent elements of $\mathcal{K}$ using
the power series notation ${\boldsymbol \varrho}^{\alpha}\equiv\exp (2\pi i{\rm Tr}({\boldsymbol \rho}\cdot \alpha ))$
(as opposed to the character notation $\phi_{\alpha}$), which has the advantage
of allowing us to interpret
Dirichlet multiplication with monomials in terms of composition: 
\[ [F]\otimes [{\boldsymbol \varrho}^{\alpha}]\;\;=\;\; [F({\boldsymbol \varrho}^{\alpha})] .\]
We denote as in the previous section $\bar{\sf N}_{\bullet}[\mathcal{K}]=\PR {\sf H}_{\bullet}[\mathcal{K}]$.

Equip ${\sf N}_{\bullet}[\mathcal{K}]$ with the Fubini-Study metric
associated to the inner-product on ${\sf H}_{\bullet}[\mathcal{K}]$.
A {\em nonlinear automorphism} 
$\Upsilon:{\sf N}_{\bullet}[\mathcal{K}]\rightarrow
{\sf N}_{\bullet}[\mathcal{K}]$
is the restriction of a graded Fubini-Study isometry of $\bar{\sf N}_{\bullet}[\mathcal{K}]$ respecting the operations $\oplus$, $\otimes$
whenever they are defined:
that is, 
\[ \Upsilon ([F]\oplus [G])\;\; =\;\; \Upsilon ([F])\oplus \Upsilon ([G]),
\quad\quad\quad\quad
\Upsilon ([F]\otimes [G])\;\; =\;\; \Upsilon ([F])\otimes \Upsilon ([G])\] and
for some permutation $\iota$ of $\mho_{\mathcal{K}}$, 
$\Upsilon (\PR{\sf H}_{\boldsymbol \vartheta}[\mathcal{K}]) = \PR{\sf H}_{\iota ({\boldsymbol \vartheta} )}[\mathcal{K}]$ for all ${\boldsymbol \vartheta} \in \mho_{\mathcal{K}}$.

For example, let $\mathcal{L}/\mathcal{K}$ be Galois.  Then the Galois group ${\sf Gal}(\mathcal{L}/\mathcal{K})$ acts on ${\sf H}_{\bullet}[\mathcal{L}]$ by: 
\[  \sum a_{\alpha}\rho^{\alpha}\longmapsto   \sum a_{\alpha}\rho^{\sigma( \alpha )}=
\sum a_{\sigma^{-1} (\alpha )}\rho^{ \alpha } \]
for $\sigma\in {\sf Gal}(\mathcal{L}/\mathcal{K})$.  This action permutes the coefficient set, hence
acts by isometries.  Viewing the action on $\mathcal{L}_{\infty}$, where it simply permutes coordinates, we see that there is an induced action on the sign set $\mho_{\mathcal{L}}$, 
so that for any ${\boldsymbol \vartheta} \in \mho_{\mathcal{L}}$, we have $\sigma (\mathcal{L}^{\boldsymbol \vartheta})=\mathcal{L}^{\sigma ({\boldsymbol \vartheta} )}$.  Thus $\sigma$
permutes the grading of ${\sf H}_{\bullet}[\mathcal{L}]$.  Finally, $\sigma$
preserves the elements of trace zero and
commutes with multiplication by elements of $\C^{\ast}$, hence induces a nonlinear automorphism
\[  \sigma :{\sf N}_{\bullet}[\mathcal{L}]\longrightarrow {\sf N}_{\bullet}[\mathcal{L}] \]
that is trivial on ${\sf N}_{\bullet}[\mathcal{K}]$.

If $\mathcal{L}/\mathcal{K}$ is a field extension of 
number fields of possibly infinite degree over $\Q$, denote by 
\[{\sf Gal}\big({\sf N}_{\bullet}[\mathcal{L}]/{\sf N}_{\bullet}[\mathcal{K}]\big) \]
the group of nonlinear automorphisms of 
${\sf N}_{\bullet}[\mathcal{L}]$ fixing
the sub nonlinear field ${\sf N}_{\bullet}[\mathcal{K}]$, and
by 
\[{\sf Gal}\big({\sf N}_{\bullet}[\mathcal{K}]\big/\mathcal{K}\big)\] the group of 
nonlinear automorphisms of ${\sf N}_{\bullet}[\mathcal{K}]$
fixing $\mathcal{K}$.  

We recall the following theorem of Wigner \cite{Wi}.
\begin{wigtheo} Let ${\bf H}$ be
a complex Hilbert space, $\PR{\bf H}=
({\bf H}\setminus\{ 0\})/\C^{\ast}$ its projectivization.
Let $[h]:\PR{\bf H}\rightarrow
\PR{\bf H}$ be a bijection preserving the Fubini-Study
metric.  Then $[h]$ is the projectivization 
of a unitary or anti-unitary linear map $h:{\bf H}\rightarrow{\bf H}$.
\end{wigtheo}

\begin{theo}\label{GalIsOne}  Let $\mathcal{K}$ be a (possibly
infinite degree) number field.  Then
\[ {\sf Gal}\big({\sf N}_{\bullet}[\mathcal{K}]/\mathcal{K}\big)
\;\;\cong\;\;\{ 1\}.\]
\end{theo}

\begin{proof}  Let 
$\sigma\in {\sf Gal}\big({\sf N}_{\bullet}[\mathcal{K}]/\mathcal{K}\big)$.
By Wigner's Theorem, $\sigma$ is the projectivization
of an (anti) unitary linear map 
\[ \tilde{\sigma}: {\sf H}_{\bullet}[\mathcal{K}]\longrightarrow
{\sf H}_{\bullet}[\mathcal{K}].\]
Since $\sigma$ fixes $\mathcal{K}$, there exist multipliers $\lambda_{\alpha}\in U(1)$
with 
\[  \tilde{\sigma}({\boldsymbol \varrho}^{\alpha})\;\; =\;\; \lambda_{\alpha}\cdot {\boldsymbol \varrho}^{\alpha}. \]
But $\sigma$ respects the Cauchy and Dirichlet products, wherein we must have
that $\lambda$ is simultaneously an additive and multiplicative character: 
\[ \lambda_{\alpha_{1}+\alpha_{2}}=\lambda_{\alpha_{1}}\lambda_{\alpha_{2}}=\lambda_{\alpha_{1}\alpha_{2}},\]
clearly possible only for $\lambda$ trivial.
\end{proof}

\begin{lemm}\label{powerlemma}  Let $\mathcal{K}$ be a (possibly infinite
degree) algebraic number field Galois over $\Q$.  Given ${\boldsymbol \vartheta} \in \mho_{\mathcal{K}}$ 
let $F\in {\sf H}_{{\boldsymbol \vartheta}}[\mathcal{K}]$ satisfy the functional equation
\begin{equation}\label{functionaleqn}
F({\boldsymbol \varrho}^{r})=(F({\boldsymbol \varrho}))^{r} 
\end{equation}  
for all $r\in\R_{+}$.
Then $F\in \mathcal{K}$ i.e. there exists $\alpha\in \mathcal{K}$ with $F({\boldsymbol \varrho})=
{\boldsymbol \varrho}^{\alpha}$
\end{lemm}

\begin{proof}  We first consider the case of $K=\Q$.  In this totally real case
we write $F$ as a function
of the parameter $\xi = \exp (2\pi i \tau)$ where $\tau = x+iy\in \HP_{\Q}$.   Assume first that
$F$ is holomorphic and non constant i.e.\ $F\in {\sf H}[\Q]$ = Hardy space of holomorphic functions.
We will show that
$F(\xi ) =\xi^{q}$ for some $q\in\Q_{+}$. 
Let us return to viewing $F$ as a holomorphic function of the half-plane variable $\tau$.  Let $\tau_{0}\in \HP_{\Q}\subset
\hat{\mathfrak{S}}_{\Q}$ be such that 
$F(\tau_{0})\not=0$.  Then there exists a complex number $\alpha$ with $\exp (2\pi i\alpha \tau_{0})=F(\tau_{0})$.
By (\ref{functionaleqn}), the functions $F(\tau)$ and $\exp (2\pi i\alpha \tau)$ agree on the ray $\R_{+}\cdot \tau_{0}$ hence by holomorphicity, they coincide.  Since $F$ is a function of $\hat{\mathfrak{S}}_{\Q}$, it follows that
$\alpha=q\in\Q_{+}$ hence $F(\xi ) =\xi^{q}$.  For $F\in {\sf H}_{-}[\Q]$ = Hardy space of anti-holomorphic functions, the argument is the same: we use the anti-holomorphic exponential
$\exp (2\pi i\alpha \bar{\tau})$ and the fact that anti-holomorphic functions agreeing on a co-dimension 1 subspace coincide.

Now let $K/\Q$ be totally real of finite degree.  Without loss of generality we may assume that
$F\in {\sf H}[K]$ = the Hardy space of holomorphic functions.  Let $\Delta\subset \HP_{K}$
be the diagonal hyperbolic sub-plane, which is
the dense leaf of the image of $\hat{\mathfrak{S}}_{\Q}$ under the diagonal embedding
$\hat{\mathfrak{S}}_{\Q}\hookrightarrow \hat{\mathfrak{S}}_{K}$. Then by the previous paragraph, the restriction of $F$
to $\Delta$ is an exponential
\[  \exp (2\pi i q\tau) \]
for $q\in \Q_{+}$.  Let $\Delta_{1}\subset \HP_{K}$ be the diagonal
\[  \Delta_{1} =\{ (\tau,\dots , \tau,\tilde{\tau}) |\; \tau,\tilde{\tau}\in \HP \}\supset \Delta.  \]
For each $\tau$ fixed, we can (using the same argument employed in
the previous paragraph) write the function $\tilde{\tau}\mapsto F(\tau,\dots , \tau,\tilde{\tau})$
as
\[  F(\tau,\dots , \tau,\tilde{\tau}) =   \exp (2\pi i q\tau) \cdot \exp (2\pi i \beta (\tau )\tilde{\tau} )\]
where $\beta (\tau )\in \C$.   Moreover we have
\[  F(r\cdot (\tau,\dots , \tau,\tilde{\tau})) =   \exp (2\pi i qr\tau) \cdot \exp (2\pi i \beta (\tau)r\tilde{\tau} )\]
by hypothesis.  As we vary $\tau$, $\beta (\tau)$ varies holomorphically
and we obtain that on a real codimension 1 
subspace of $\Delta_{1}$,
\[ F(\tau,\dots , \tau,\tilde{\tau}) = F_{1}(\tau,\dots ,\tau,\tilde{\tau}) := \exp (2\pi i q\tau) \cdot \exp (2\pi i \beta (\tau)\tilde{\tau}) \] 
hence they are equal on $\Delta_{1}$.  But this
means that
$F_{1}$ must also satisfy the functional equation $F_{1}(r\cdot (\tau,\dots , \tau,\tilde{\tau})) = 
F_{1}(\tau,\dots ,\tau,\tilde{\tau})^{r}$, which
implies that $\beta  (\tau)=\beta$ is a constant.  Inductively, we obtain
that $F$ restricted to $\HP_{K}$ is an exponential function, and since $F$ is a function of $\hat{\mathfrak{S}}_{K}$, this restriction must be of the form
${\boldsymbol \xi}^{\alpha}=\exp (2\pi i {\rm Tr}(\alpha\cdot \boldsymbol{\tau}))$ for $\alpha\in K_{+}$.
The case of a  totally complex field extension, one of mixed type, or one of infinite degree, 
is handled in exactly the same manner.

\end{proof}
 
\begin{theo}\label{GalIsGal} Let $\mathcal{L}/\mathcal{K}$ be
a Galois extension of (possibly infinite degree) algebraic number fields. Then 
\[ {\sf Gal}\big({\sf N}_{\bullet}[\mathcal{L}]/
{\sf N}_{\bullet}[\mathcal{K}]\big) 
\;\;\cong\;\; {\sf Gal}(\mathcal{L}/\mathcal{K} ).\] 
\end{theo}

\begin{proof}  Let $\sigma\in 
{\sf Gal}({\sf N}_{\bullet}[\mathcal{L}]/{\sf N}_{\bullet}[\mathcal{K}])$.
We begin by showing
that $\sigma (\mathcal{L})=\mathcal{L}$, 
where $\mathcal{L}$ is identified with the field of monomials
$[{\boldsymbol \varrho}^{\alpha}]$, $\alpha\in\mathcal{L}$.  Note first that we have
already $\sigma (\mathcal{K})=\mathcal{K}$.
Since $\sigma (\mathcal{L})$ is a field, all of its elements obey
the distributive law.  Thus, given any $[F]\in\sigma (\mathcal{L})$, since 
$\sigma( [{\boldsymbol \varrho}^{m}])=[{\boldsymbol \varrho}^{m}]\in\mathcal{K}$, $m\in\N$, we have
\begin{eqnarray*}[F({\boldsymbol \varrho}^{m})]\;\;=\;\; 
[F]\otimes [{\boldsymbol \varrho}^{m}] & = & 
[F]\otimes \big([{\boldsymbol \varrho}]\oplus\cdots \oplus [{\boldsymbol \varrho}]\big) \\
& = &
\big( [F]\otimes [{\boldsymbol \varrho}] \big)\oplus \cdots 
\oplus \big([F]\otimes [{\boldsymbol \varrho}] \big) \\ 
& = & [F({\boldsymbol \varrho})]\oplus \cdots \oplus [F({\boldsymbol \varrho})] \\
& = & [F]^{m} ,
\end{eqnarray*}
where $[F]^{m} $ denotes the Cauchy $m$th power of $[F]$.
In fact, the same argument shows that for any $m/n\in\Q_{+}$
\[ ([F({\boldsymbol \varrho}^{m/n})])^{n}\;\;=\;\; [F]^{m}.\]
We may thus find $F\in [F]$ 
satisfying the functional
equation
$F({\boldsymbol \varrho}^{q}) = (F({\boldsymbol \varrho}))^{q} $
for all $q\in\Q_{+}$.  By continuity, this
extends to $\R_{+}$. 
Note that since $\sigma$ respects
the grading and each element of $\mathcal{L}$ is homogeneous
(is contained in a fixed projective summand)
then $[F]$ is also homogeneous.  Thus,
by Lemma~\ref{powerlemma}, it follows that $F\in \mathcal{L}$ and
$\sigma (\mathcal{L})=\mathcal{L}$.   We induce in this way a homomorphism
\[ \Pi : 
{\sf Gal}\big({\sf N}_{\bullet}[\mathcal{L}]/{\sf N}_{\bullet}[\mathcal{K}]\big) 
\;\;\longrightarrow\;\; {\sf Gal}(\mathcal{L}/\mathcal{K} ),
\quad\quad  \Pi (\sigma )=\sigma|_{\mathcal{L}} .\]
Note that $\Pi$ is clearly onto, as we have already observed that any $\sigma\in {\sf Gal}(\mathcal{L}/\mathcal{K} )$
generates an automorphism of ${\sf N}_{\bullet}[\mathcal{L}]$ fixing
${\sf N}_{\bullet}[\mathcal{K}]$ via ${\boldsymbol \varrho}^{\alpha}\mapsto {\boldsymbol \varrho}^{\sigma (\alpha )}$.
Suppose now that $\Pi (\sigma )=1$ for $\sigma$
in ${\sf Gal}\big({\sf N}_{\bullet}[\mathcal{L}]/{\sf N}_{\bullet}[\mathcal{K}]\big) $.  
Then $\sigma\in 
{\sf Gal}\big({\sf N}_{\bullet}[\mathcal{L}]/\mathcal{L} \big)$, but by Theorem \ref{GalIsOne},
the latter group is trivial.
\end{proof}

We now concentrate on the case of a finite Galois extension $K/\Q$ and
consider each of the operations
$\oplus$ and $\otimes$ separately.  
We will work with the nonlinear number field $\bar{\sf N}_{\bullet}[K]$ = 
$\PR{\sf H}_{\bullet}[K]$.
Let ${\sf Gal}_{\oplus}\big(\bar{\sf N}_{\bullet}[K]/K\big)$
denote those isometries fixing $K$ and homomorphic with respect
to $\oplus$ only. 
${\sf Gal}_{\otimes}\big(\bar{\sf N}_{\bullet}[K]/K\big)$
is defined similarly.

Denote by ${\sf U}\big({\sf H}_{\bullet}[K]\big)$
the group of unitary operators of ${\sf H}_{\bullet}[K]$.
The action of ${\bf r}\in K_{\infty}$ by 
translation in $\HP_{K}$, 
${\bf z}\mapsto {\bf z}+{\bf r}$, induces an action
on ${\sf H}_{\bullet}[K]$ by
\[  \Phi_{\bf r}(F)\;\; =\;\;   
\sum a_{\alpha}\exp (2\pi i{\sf Tr}( \alpha {\bf r}) ){\boldsymbol \varrho}^{\alpha} 
\]
for $F=\sum_{\alpha}a_{\alpha}{\boldsymbol \varrho}^{\alpha}$, yielding a faithful representation
\[ \Phi :K_{\infty}\;\longrightarrow\;  
{\sf U}\big({\sf H}_{\bullet}[K]\big) .\]

\begin{prop}\label{projflow}
The projectivization $[\Phi ]$ of $\Phi$
defines a monomorphism 
\[ [\Phi ]: K_{\infty}\;\hookrightarrow \;
{\sf Gal}_{\oplus}\big(\bar{\sf N}_{\bullet}[K]/K\big)
. \]
\end{prop}

\begin{proof}  For $[{\boldsymbol \varrho}^{\alpha}]\in K$ and ${\bf r}\in K_{\infty}$, 
$[\Phi ]_{\bf r}([{\boldsymbol \varrho}^{\alpha}])=
[\exp (2\pi i{\rm Tr}( \alpha {\bf r})  ){\boldsymbol \varrho}^{\alpha}]=[{\boldsymbol \varrho}^{\alpha}]$.
For any $[F]$, $[G]\in {\sf N}_{\bullet}[K]$, let $[f]$, $[g]$
be the projective classes of their boundary functions.  Then
$[\Phi ]_{\bf r}([F]\oplus [G])$ is the element
of ${\sf N}_{\bullet}[K]$ whose boundary function is
\begin{eqnarray*}
[\Phi ]_{\bf r}([f]\oplus [g]) & = & 
\left[ \sum_{\alpha}\Big(\sum_{\alpha_{1}+\alpha_{2}=\alpha }a_{\alpha_{1}}b_{\alpha_{2}}\Big)
\exp (2\pi i{\sf Tr}( \alpha {\bf r}) ){\boldsymbol \zeta}^{\alpha} \right]\\
& \\
& = & 
\left[ \sum_{\alpha}\Big(\sum_{\alpha_{1}+\alpha_{2}=\alpha }
a_{\alpha_{1}}\exp (2\pi i{\rm Tr}( \alpha_{1} {\bf r}))b_{\alpha_{2}}
\exp (2\pi i{\sf Tr}( \alpha_{2} {\bf r} ))\Big)
{\boldsymbol \zeta}^{\alpha}\right] \\
& \\
& = & [\Phi ]_{\bf r}([f])\oplus [\Phi ]_{\bf r}([g]) ,
\end{eqnarray*}
which is the boundary function of $[\Phi ]_{\bf r}([F])\oplus [\Phi ]_{\bf r}([G])$.
\end{proof}

In like fashion, we may define a flow on
$\bar{\sf N}_{\bullet}[K]$ respecting $\otimes$ as follows.
For a vector
${\bf z}\in K_{\infty}$ we denote by $\log |{\bf z}|$
the vector 
\[  \big(\log |z_{\nu_{1}}|,\dots ,\log |z_{\nu_{d}}|\big)\]
when $K$ is real, or when $K$ is complex
\[  \big(\log |z_{\mu_{1}}|,\log |z_{\mu_{1}}|\dots ,\log |z_{\mu_{s}}|,\log |z_{\mu_{s}}|\big).\]
Then for
$F=\sum_{\alpha}a_{\alpha}{\boldsymbol \varrho}^{\alpha}$ we define
\[  \Psi_{\bf r} (F)\;\; =\;\; 
\sum_{\alpha\in K}a_{\alpha}
\exp \big( 2\pi i{\sf Tr}({\bf r} \log | \alpha | )\big) 
{\boldsymbol \varrho}^{\alpha}.\]
This defines a faithful representation
\[ \Psi :K_{\infty}\;\longrightarrow\;  
{\sf U}\big({\sf H}_{\bullet}[K]\big) .\]
The following is proved exactly as Proposition~\ref{projflow}:

\begin{prop}\label{projflow2}  The projectivization $[\Psi ]$
of $\Psi$ defines a monomorphism
\[ [\Psi ]:  K_{\infty}\;\hookrightarrow \;
{\sf Gal}_{\otimes}\big(\bar{\sf N}_{\bullet}[K]/K\big)
. \]
\end{prop}

Recall that the idele class group of $\Q$, $C_{\Q}$, may be
identified with $\R^{\ast}_{+}\times {\rm Gal}(\bar{\Q}^{\rm ab}/\Q)$.  We have
the following Corollary to Propositions \ref{projflow} and \ref{projflow2}:

\begin{coro} There are monomorphisms
\[ {\sf C}_{\Q}\;\hookrightarrow\; 
{\sf Gal}_{\oplus}\big(\bar{\sf N}_{\bullet}[\Q^{\rm ab}]/\Q\big) ,
\quad\quad\quad\quad
{\sf C}_{\Q }\;\hookrightarrow\; 
{\sf Gal}_{\otimes}\big(\bar{\sf N}_{\bullet}[\Q^{\rm ab}]/\Q\big) .
\]
\end{coro}

We end by noting that the above structures have an interpretation within the
von Neumann depiction of quantum mechanics.
The space $\bar{\sf N}_{\bullet}[K]$
may be viewed as the space of states of a quantum mechanical system
with $d$ degrees of freedom.
We view each of the flows $\Phi$ and $\Psi$ as generating
the coordinate observables of two distinct physical systems and write
\[ \Phi_{\bf r}\;\;=\;\; \exp (2\pi i\langle H_{\oplus},{\bf r}\rangle ) ,\quad\quad\quad\quad 
  \Psi_{\bf r}\;\;=\;\; \exp (2\pi i \langle H_{\otimes}, {\bf r}\rangle),\]
where $H_{\oplus}=\textrm{diag}[ q ]_{q\in K}$
and $H_{\otimes}=\textrm{diag}[ \log |q| ]_{q\in K}$ are the associated
(vector-valued)
Hamiltonian operators. The set of stationary states for 
each of $H_{\oplus}$ and $H_{\otimes}$ is the field $K$.

Each $[F]\in \bar{\sf N}_{\bullet}[K]$ defines a Cauchy multiplication
operator $M_{\oplus}([F])$ and a Dirichlet multiplication operator
$M_{\otimes}([F])$. It is easy to see that if $[F]$ has a representative
$F$ whose Fourier coefficients are real, then $M_{\oplus}([F])$ defines
an observable {\em i.e.}\
the projectivization of a self-adjoint operator of ${\sf H}_{\bullet}[K]$.  
This is not true of $M_{\otimes}([F])$ due
to the fact that the Haar measure on $\SO_{K}$ --
which we use to define the Hardy inner-product -- is invariant with respect
to addition but not multiplication.  This suggests that regarding
the system defined by $H_{\otimes}$,
it may be more
natural to use the Hardy space of functions holomorphic
on a hyperbolized idele class group ${\sf C}_{K}$
(with its multiplicatively invariant Haar measure).
A formal computation shows that the operators
$M_{\otimes}([F])$ are self-adjoint for the ``idelic'' Hardy inner-product
when the Fourier coeficients of some representative $F$
are real.  We suspect that 
the eigenvalues of $M_{\otimes}([F])$
are equal or related in a straightforward manner 
to the imaginary parts of the zeros of a meromorphic
extension of a Dirichlet type series corresponding
to $[F]$.

\section{Appendix: Errata to the Published Version \cite{GeVeP}}

In what follows, the revised version (that is, {\it this version}) is denoted [R].

\vspace{3mm}

\noindent {\bf 1.}
The definition given of nonlinear number field was erroneous.  
On page 582, line 15 of \cite{GeVeP} the definition of $\C^{\ast} [K]$ should read:
\begin{equation}\label{corr}  \C^{\ast} [K] = \C [K] -I_{K}
\end{equation}
and not ``$\C [K]/I_{K}$''.  
This correction must be carried out in the more general
definition {\it e.g.}\ on page 586 line 30 of \cite{GeVeP} it should read
$ {\sf H}^{\ast}_{\bullet} [K] =  {\sf H}_{\bullet}[K] -{\sf I}_{K} $, where ${\sf I}_{K}$ is the kernel
of the trace map on its domain of definition.  See \S 9 of [R].

\vspace{3mm}

\noindent {\bf 2.} The discussion of the Dirichlet product structure of 
${\sf Char}(\SO_{K})$ (part 1 of Theorem 5, page 581 of \cite{GeVeP}) was incomplete.   See \S 7 of  [R].

 \vspace{3mm}

\noindent {\bf 3.}    In part 2 of Theorem 5 of \cite{GeVeP},  the additive group
$(O_{K},+)$ was mistakenly identified with the character group 
${\sf Char}(\T_{K})$ of
the Minkowski torus $\T_{K}$.
Rather, it is the inverse different $\mathfrak{d}_{K}^{-1}\supset O_{K}$ which is identified
with ${\sf Char}(\T_{K})$.  The character group ${\sf Char}(\T_{K})$ is thus an $O_{K}$-module
extension of $O_{K}$.  See See \S 7 of  [R].

 \vspace{3mm}

\noindent {\bf 4.}  On line 15, page 583 of \cite{GeVeP}, it was incorrectly asserted that the Cauchy and Dirichlet products are fully defined on the Hilbert space $L^{2}(\SO_{K},\C )$.
These operations only extend partially so that for each element $f\in 
L^{2}(\SO_{K},\C )$ one must specify the Cauchy and Dirichlet domains
${\sf Dom}_{\oplus}(f)$, ${\sf Dom}_{\otimes}(f)$ of elements $g$ for which
$f\oplus g$ resp. $f\otimes g$ make sense.    The definition
of nonlinear number field, which appears in Definition 2 on page
587 of \cite{GeVeP}, must be adjusted accordingly by replacing everywhere the phrase ``double semigroup''
by ``partial double semigroup'' to take into account this correction.  See \S\S 8,9 of [R].

 \vspace{3mm}

\noindent {\bf 5.}  The proof of Lemma 1 on page 589 of \cite{GeVeP} was incorrect.  A correct proof can be
found in \S 10 of [R] (where it is known as ``Lemma 2'').

\vspace{3mm}

\noindent {\bf 6.}  The discussion of the idele class group found
on pages 591-592 of  \cite{GeVeP} is only valid for $K=\Q$, so that the idele class group
${\sf C}_{K}$  should be replaced by ${\sf C}_{\Q}$.

\vspace{3mm}

\noindent {\bf 7.} Apart from implementing the above corrections, \S 6 of the revised version [R] 
contains an enhancement of the hyperbolization $\HP_{K}$ (page 579-580 of  \cite{GeVeP}).   Furthermore, the $\theta$-holomorphic grading of functions on $\HP_{K}$
(described on page 585 of  \cite{GeVeP}) as been expanded 
by a  {\it complex sign set} along the complex places.  
This is described 
in \S 8.2 of [R].

\bibliographystyle{amsalpha}

\begin{thebibliography}{A}



\bibitem [1]{Ca-Co} A. Candel \& L. Conlon, {\em Foliations I}, 
American Mathematical Society, 2000.

\bibitem [2]{Ca-Fr} J.W.S. Cassels \& A. Fr\"{o}hlich (ed.), 
{\it Algebraic Number Theory}, London Mathematical Society; 2nd Revised edition edition,
London, UK, 2010. 

\bibitem [3]{Co} A. Connes, Noncommutative geometry and the Riemann zeta function, 
in {\em Mathematics: frontiers and perspectives}, 
ed. V. Arnold, M. Atiyah,
P. Lax \& B. Mazur, 
American Mathematical Society, 2000, pp. 35--54. 

\bibitem [4p]{GeVeP} Gendron, T.M. \& Verjovksy, A., {\it Geometric Galois
theory, nonlinear number fields and a Galois group interpretation of the idele
class group.} {\it Internat. J. Math.}  {\bf 16}  (2005),  no. 6, 567--593.


\bibitem [4e]{GeVeE} Gendron, T.M. \& Verjovksy, A., Errata to  {\it Geometric Galois
theory, nonlinear number fields and a Galois group interpretation of the idele
class group.}  to appear in  {\it Internat. J. Math.}

\bibitem [5]{GeVeN} Gendron, T.M. \& Verjovksy, A., Notes on nonlinear number fields.
arXiv:1007.3070v1  [math.NT].

\bibitem [6]{Gh} E. Ghys,
Laminations par surfaces de Riemann, in
{\it Dynamique et g\'{e}om\'{e}trie complexes},   
Panor. Synth\`{e}ses, {\bf 8}, 
Soc. Math. France, Paris, 1999, pp. 49--95. 


\bibitem [7]{Go} C. Godbillon,{\it Feuilletages. 
\'{E}tudes g\'{e}om\'{e}triques}, Birkh\"{a}user, 1991.

\bibitem [8]{Ko} H. Koch, {\it Algebraic Number Theory},
Encyclopaedia of Mathematical Sciences {\bf 62}, Number Theory II,
Springer-Verlag, 1992.


\bibitem [9]{La} S. Lang, {\it Algebra}, Graduate Texts in Mathematics {\bf 73},
Springer-Verlag, 1992.


\bibitem [10]{Mo-Sc} C.C. Moore \& C. Schochet, 
{\it Global Analysis on Foliated Spaces}, 
Mathematical Sciences Research Institute Publications {\bf 9}, 
Springer-Verlag, 1988.

\bibitem [11]{Ne} Neukirch, {\it Algebraic Number Theory},
Grundlehren der Mathematischen Wissenschaften {\bf 322}. Springer-Verlag, Berlin, 1999.

\bibitem [12]{Pr-Se} A. Presley \& G. Segal,
{\it Loop Groups}, Oxford University Press, 1986. 

\bibitem [13]{Ra-Va}  D. Ramakrishnan \& R.J. Valenza, Graduate Texts in Mathematics {\bf 186},
{\it Fourier Analysis on Number Fields}, 
Springer-Verlag, 1999.

\bibitem [14]{Ta} Tate, J. T.,
Fourier analysis in number fields, and Hecke's zeta-functions. in {\it Algebraic Number Theory (Proc. Instructional Conf., Brighton, 1965)} (ed. Cassels \& Fr\"{o}hlich), pp. 305--347, London Mathematical Society; 2nd Revised edition, London, UK, 2010.

\bibitem [15]{We1} A. Weil, On the Riemann hypothesis in function fields, 
{\it Proc. Nat. Acad. Sci. U. S. A.} {\bf 27} (1941), 345--347.

\bibitem [16]{We2} A. Weil, Sur la th\'{e}orie de corps de classes,
\textit{J. Math. Soc. Japan} {\bf 3} (1951).

\bibitem [17]{Wi} E. Wigner, {\it Gruppentheorie und ihre Anwendung auf die 
Quantenmechanik der Atomspektren},  
J. W. Edwards, 1944. 

\end{thebibliography}

\end{document}